\documentclass{article}

\usepackage{PRIMEarxiv}
\usepackage{graphicx}
\usepackage{subcaption}
\usepackage{color}
\usepackage{float} % fix table inside context 
\usepackage{fullpage}
\usepackage{adjustbox}
\usepackage{amsfonts}
\usepackage{amsmath}
\usepackage{mathtools}
\usepackage{accents}
\usepackage{amssymb}
\usepackage{amsbsy}
\usepackage{amsthm}
\usepackage{graphicx}
\usepackage{colordvi}
\usepackage{graphics}
\usepackage[dvipsnames]{xcolor}
\usepackage{bm}
\usepackage{algorithm}
\usepackage{algorithmic}
\usepackage{verbatim} 
\usepackage{braket}

\usepackage{siunitx}

\usepackage{amsfonts}

\usepackage{dsfont}
\usepackage{graphicx}
\usepackage{algorithmic}
\usepackage{subcaption}
\usepackage{accents}

\newcommand{\Int}[4]{\displaystyle{\int\limits_{#1}^{#2} {#3} d{#4}}}

\newcommand{\nn}{\boldsymbol{n}}
\newcommand{\XTO}{X_{\mathcal{T}_+}}
\newcommand{\XTT}{X_{\mathcal{T}_-}}

\newcommand{\specRad}{\rho(K_+^{-1}SK_-^{-1}D)}
\newcommand{\Tpl}{\boldsymbol{T}_+^h}
\newcommand{\Tmi}{\boldsymbol{T}_-^h}
\newcommand{\fp}{\boldsymbol{f}_+}
\newcommand{\fm}{\boldsymbol{f}_-}
\newcommand{\tminus}{T_-^{k+1/2}}
\newcommand{\tplusw}{\widetilde{T}_+^{k+1}}
\newcommand{\tplus}{T_+^{k+1}}

\newcommand{\tminusstar}{T_{-}^{*}}
\newcommand{\tpluswstar}{\widehat{T}_{+}^{*}}
\newcommand{\tplusstar}{T_{+}^{*}}
\newcommand{\tlaststar}{T_{+}^{*}}

\newcommand{\wlast}{w_+^k}
\newcommand{\wminus}{w_-^{k+1/2}}
\newcommand{\wplus}{\widetilde{w}_+^{k+1}}

\newcommand{\GmD}{\Gamma_{-\,D}}

\newcommand{\Edgesi}{\mathcal{E}_i}

\newcommand{\duality}[3]{\silvia{\left\langle #2,#3 \right\rangle_{#1}}}	

\newtheorem{prop}{Proposition}

\newtheorem{problem}{Problem}

\newtheorem{theorem}{Theorem}[section]
\newtheorem{thm}{Theorem}[section]

\newtheorem{remark}[thm]{Remark}

\def\alex#1{
 {\color{black}#1}
 }

\def\silvia#1{
{\color{black}#1}
}

\title{A theoretical and numerical analysis of a Dirichlet-Neumann domain decomposition method for diffusion problems in heterogeneous media
}

\author{
  Alex Viguerie \\
  Division of Mathematics \\
  Gran Sasso Science Institute \\
  L`Aquila, AQ 67100, Italy\\
  \texttt{alexander.viguerie@gssi.it} \\
   \And
  Silvia Bertoluzza \\
  CNR Imati Enrico Magenes \\
  Pavia, PV 27100, Italy  \\
  \texttt{silva.bertoluzza@imati.cnr.it} \\
   \AND
  Alessandro Veneziani \\
  Department of Mathematics\\
  Department of Computer Science \\
  Emory University \\
  Atlanta, GA 30322, USA \\
  \And
  Ferdinando Auricchio \\
  Dipartimento di Ingegneria Civile e Architettura \\
  Universita degli Studi di Pavia \\
 Pavia, PV 27100, Italy
}

\begin{document}
\maketitle

\begin{abstract}
Problems with localized nonhomogeneous material properties present well-known challenges for numerical simulations. In particular, such problems may feature large differences in length scales, causing difficulties with meshing and preconditioning. These difficulties are increased if the region of localized dynamics changes in time. Overlapping domain decomposition methods, which split the problem at the continuous level, show promise due to their ease of implementation and computational efficiency. Accordingly, the present work aims to further develop the mathematical theory of such methods at both the continuous and discrete levels. For the continuous formulation of the problem, we provide a full %well-posedness and 
convergence analysis. For the discrete problem, we show how the described method may be interpreted as a Gauss-Seidel scheme or as a Neumann series approximation, establishing a convergence criterion in terms of the spectral radius of the system. We then provide a spectral scaling argument and provide numerical evidence for its justification.

% Please include a maximum of seven keywords
%\keywords{Two-level DD methods, finite element methods, domain decomposition, thermal problems}
\end{abstract}

\section{Introduction}\label{sec:intro}
Many problems in science and engineering feature localized dynamics in which material properties may differ across various regions of the problem domain. Such problems are common, for example, in additive manufacturing (AM) \cite{ BRDP2018, GMWP2012, HYDY2016,IM2016, KAFHKKR2015, KOCZR2018, LRMR2016, PPRZMHBS2015, PPRZMHBS2016, RSM2014, RSZ2018, TT2017,WLGNMR2017}, where accurate simulation of the problem requires adequate resolution of phenomena at both the micrometer and millimeter scale. Problems of this type are a well-known source of numerical difficulty, as the presence of multiple length scales, irregular geometries, and \silvia{the presence of interfaces arising from abrupt changes in material properties}
%large gradients arising from changes in material properties
often complicates the meshing and simulation process.  \alex{Such problems are compounded if the regions requiring small-scale resolution change in time \cite{VA2019}}.
\par To address these challenges, \alex{many methods have been proposed}. The use of classical methods involves remeshing at each time step, which requires difficult-to-implement and expensive refinement-and-derefinement algorithms \cite{Patil2015, carraturo2019a, kollmannsberger_hierarchical_2018, LI2019100903, baiges2020adaptive}.  \silvia{Expensive meshing/remeshing can be avoided by resorting to unfitted methods, such as GFEM/XFEM \cite{Kergrene:Babauska:Banerjee:2016, Belytschko:Gracie:Ventura:2009} (and, more generally, partition of unit methods, (\cite{Melenk:Babuska:1996}),} \textcolor{black}{Immersed Finite Element \cite{peskin1977numerical, liu2007mathematical},} \silvia{CutFEM (\cite{Burman:Elferveson:Hansbo:etal:2019}), or hierarchical methods such as the finite cell method with local enrichment \cite{Joulaian:Meysam:Duester:2013} or} $hp$-$d$ methods \cite{ DNR2007, PP2009, R1992, SDR2012, ZBEFKR2016}. Alternatively, one can resort to methods that facilitate the remeshing procedure by allowing polygonal elements, with possibly curved edges/faces, such as the virtual element method \cite{de2016nonconforming, beirao2014hitchhiker, brezzi2005family, Long:Huayi:Wen:2017, Aldakheel:Fadu:Hudobivnik:etal:2009}.

%These include the \textit{fat boundary} \cite{BMM2005, BIM2011, Maury2001},   and 

\silvia{Combining a domain decomposition method with a fictitious domain approach, results in the method proposed in \cite{VBA2019}, referred to herein as the Two-level DD method. More precisely, already at the continuous level, the domain is split as the union of two regions, corresponding to two different materials, each one homogeneous. The problem is then tackled by a non-overlapping Dirichlet-Neumann DD method. Assuming that one of the two regions is small, the corresponding problem (which we refer to as the \textit{local problem}) is solved by a standard finite element method on a fitted mesh. The problem in the larger region (which we refer to as the \textit{global problem}) is solved by resorting to a fictitious domain approach: both coefficients and right hand side are suitably extended to the whole domain, and, in the spirit of the Fat Boundary method \cite{Maury2001, BIM2011},  information on the jump of the normal flux on the interior interface is retrieved from the solution of the local problem and injected as a data in the global problem.  The resulting algorithm exhibits many of the desirable characteristics one may expect from a domain decomposition approach, as it is easy to implement \alex{and precondition}, and allows for the  employment of quasi uniform meshes, of possibly different size, for both local and global problems. \alex{The regularity of both meshes and unchanging topology of the global domain, usually corresponding to large scales, make this approach particularly well-suited for problems in which the (small) region, where local-scale phenomena occur, evolves in time, as in \cite{VA2019}. In such cases, contrary to most of the other approaches considered in the literature for the solution of problem of this kind, the approach of \cite{VA2019, VBA2019} allows to completely avoid not only remeshing, but also recomputing the entries of the local and global stiffness matrices, 
%, representing an alternative to difficult-to-implement and expensive refinement-and-derefinement algorithms \cite{Patil2015, carraturo2019a, kollmannsberger_hierarchical_2018, LI2019100903, baiges2020adaptive}
as changes in the domain configurations can be handled by comparatively simple mesh translations for the local domain.}} 
%These aspects were further demonstrated in \cite{VA2019}.
\par \alex{Both domain decomposition (DD) and fictitious domain/immersed boundary methods have a long history and are a well-studied topic in the literature \cite{lions1988schwarz, toselli2004domain, lions1990schwarz, QV1999, gunzburger1999optimization, Glowinski:Pan:Periaux:1994, Glowinski:Pan:etal:1999, WOS:A1994PB02700004}. Indeed, the aforementioned Fat-boundary and $hp$-$d$ methods can themselves be regarded as variants within the classical DD framework \cite{R1992, BMG2012}. Domain decomposition methods have many variants, and are generally classified according to whether the subdomains are overlapping or not, as well as on the way in which information is transferred between the subdomains. For non-overlapping DD, the information exchange is carried out through suitable boundary conditions at the interface
such ad Dirichlet-Neumann or Robin-Robin, the latter being at the basis of the optimized Schwarz method. Information between subdomains can also be exchanged by the use of Lagrange multiplier, as in the Mortar method, or by introducing an auxiliary variable (which might be regarded as a control variable) on the interface \cite{QV1999, cote2005comparison, gunzburger1999optimization, Bertoluzza:Brezzi:Sangalli:2007}. Approaches of this type are natural for problems exhibiting particular geometric characteristics, and can be used both directly as numerical solvers or as preconditioners \cite{giraud2006algebraic, QV1999, da2006positive, da2010robust, da2012overlapping,da2013isogeometric}. }

\par \alex{Though domain decomposition algorithms have indeed been used for problems of the type studied in the current work, combining them with fictitious domain type methods to tackle problems in which the domains are separated by differences in physical materials requires some care in the design of the transmission conditions, leading to the formulations expressed in the current work. Though the application of such an approach in \cite{VA2019, VBA2019} shows potential from a numerical point of view, significant theoretical questions remain. While the Two-level DD formulation was shown in \cite{VBA2019} to be consistent with the original problem formulation, the convergence behavior of the Schwarz type iterations  was not proven. At the continuous level, it is known that the differences in the material properties may affect the convergence of a DD method, requiring the introduction of under-relaxation or similar regularization techniques. At the discrete level, the potentially large difference in mesh resolutions and/or discretization methods necessary across the different domains may also cause convergence issues. In the present work, we begin to address the aforementioned issues, starting with diffusion problems in heterogeneous media with piece-wise constant coefficients, for which we provide a theoretical foundation for the practical application of Dirichlet-Neumann non-overlapping domain decomposition, coupled with a fictitious domain method for the global problem. We seek to clearly prove and establish appropriate convergence behavior and conditions analytically. Numerically, we aim to investigate and better understand the effect of different materials and mesh resolutions on such convergence behavior.}

\par  \alex{The article is outlined as follows.} We first introduce the model problem and the corresponding Two-level DD formulation shown in \cite{VBA2019} (Sect. \ref{sec:statement}). We then proceed to provide a full
 %well-posedness and
  convergence analysis (Sect. \ref{sec:conv}) . Next, we shift our attention to the discrete problem.  We will show that the discrete Two-level DD method can be equivalently interpreted as a Gauss-Seidel method or a truncated Neumann series, implying a spectral convergence condition (Sect. \ref{sec:discrete}). Accordingly, we then provide a scaling argument for the system eigenvalues in terms of the problem parameters, which is validated through a series of numerical simulations on two- and three-dimensional problems (Sect. \ref{sec:numexp}).
  Follow-ups of the present work are drawn in Sect. \ref{sec:concl}.

\begin{figure}
	\centering
	\includegraphics[width=.7\linewidth]{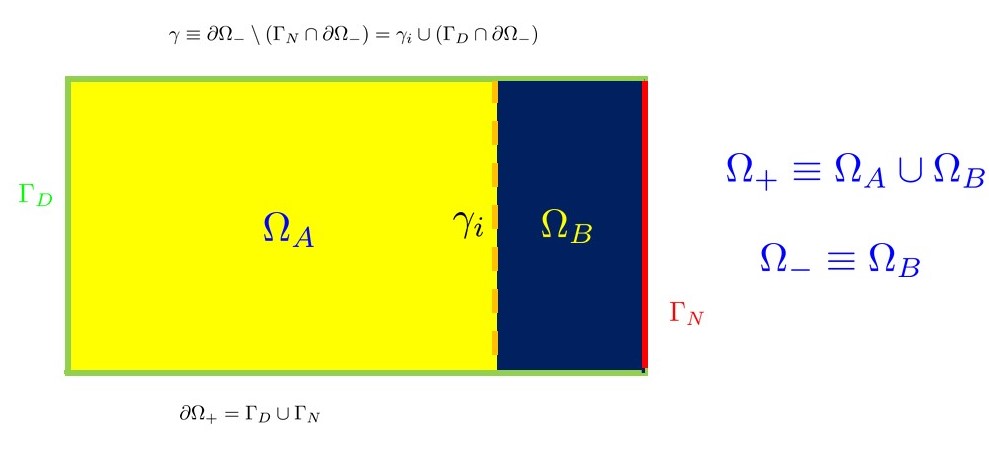}
	\caption{Domain and terminology. In the colored version, $\Gamma_N$ is denoted in red and $\Gamma_D$ in green.}
	\label{fig:domain}
\end{figure}

\section{Statement of the Problem}
\label{sec:statement}
We consider a domain given by $\Omega_+ = \Omega_A \cup \Omega_B $ where $\Omega_A \cap \Omega_B = \emptyset$ (Fig. \ref{fig:domain} ). Given $f \in L^2 (\Omega_+), q \in H^{-1/2}(\Gamma_N)$ and $T_D \in H^{1/2}(\Gamma_D)$,
we solve the following variational problem 
on $\Omega_+$. 
\begin{problem}\label{origProblem}
{\bf(Original Problem)} Find $T \in H^1(\Omega_+)$ such that:
	\begin{equation}\label{eq:origProblem}
	\Int{\Omega_+}{}{\kappa \nabla T \cdot \nabla v}{\Omega} = \Int{\Omega_+}{}{fv}{\Omega} + \Int{\Gamma_N}{}{qv}{\Gamma}, 
		\quad \forall v \in H^1_{\Gamma_D} (\Omega_+),	
	\quad T(\Gamma_D) = T_D, 
	\end{equation}
	with
	$$
	\kappa = \left\{ \begin{array}{lc} \kappa_A & \mathrm{in} \ \Omega_A \\[5pt] \kappa_B & \mathrm{in} \ \Omega_B, \end{array} \right.
	$$
\noindent \silvia{$\kappa_A, \kappa_B$ positive constants,} where $\Gamma_D \subseteq \partial \Omega_+$ (respectively $\Gamma^N\subset \partial \Omega_+$) denote the Dirichlet (resp. Neumann) boundary, and  $H^1_{\Gamma_D}(\Omega_+) = \{ u \in H^1(\Omega_+): \ u=0 \text{ on } \Gamma_D\}$. Precisely, $T \in {\cal L} + H^1_{\Gamma_D}(\Omega_+)$, where ${\cal L}$ is an arbitrary lifting of $T_D$.
\end{problem}
\

\silvia{
In the following we let $\gamma_i = \partial\Omega_+ \setminus \partial \Omega_B$ denote the portion of the boundary of $\Omega_B$ interior to $\Omega_+$, and we assume that $\bar \gamma_i \cap \bar \Gamma_N = \emptyset$. For $\gamma \subseteq \partial\Omega_B$ we let $\duality{\gamma}{\cdot}{\cdot}$ denote the duality relation between $(H^1_{\Gamma_D}(\Omega_B)|_{\gamma})'$ and $(H^1_{\Gamma_D}(\Omega_B)|_{\gamma})$. We recall that for $T \in H^1(\Omega_B)$ with $\Delta T \in L^2(\Omega_B)$ we have that $\nabla T$ has a normal trace on $\gamma_i$ which belongs to  $(H^1_{\Gamma_D}(\Omega_B)|_{\gamma_i})'$, and that we have the following ``integration by part'' identity (which actually is, for a generic $T$ with no extra smoothness, the definition of the trace of the normal derivative): for all $v \in H^1_{\Gamma_D}(\Omega_+)$
\begin{equation}\label{intbyparts}
\duality{\gamma_i}{\nabla T \cdot \nn}{v} = \int_{\Omega_B} \nabla T \cdot \nabla v + \int_{\Omega_B} \Delta T \, v.
\end{equation}
}
Following \cite{VBA2019}, we then  split the Problem (\ref{origProblem}) into two coupled problem for the unknowns $T_+ = T$ and $T_- =T|_{\Omega_-}$: a \textit{global problem} defined on all of $\Omega_+$ and a \textit{local problem} defined on $\Omega_- = \Omega_B$ as follows:

\begin{problem}[Global/local Reformulation]\label{CoupledProblem}	Find $T_+ \in H^1(\Omega_+)$ and $T_- \in H^1(\Omega_-)$ such that
	\begin{itemize}
		\item[$\rightarrow$] \textsc{Global}
	\begin{eqnarray}\label{globalPbm}
& \Int{\Omega_+}{}{ \kappa_+ \nabla T_+ \cdot \nabla v}{\Omega}
-   (\kappa_+-\kappa_-) \duality{\gamma_i}{\nabla T_- \cdot \nn}{ v }
= \Int{\Omega_+}{}{\widetilde{f} v}{\Omega} + \Int{\Gamma_N}{}{\widetilde q v}{\Gamma}, 	\quad \forall v \in H^1_{\Gamma_D} (\Omega_+), \\
& T(\Gamma_D) = T_D 
\end{eqnarray}
with $\kappa_+ = \kappa_A$ on the whole $\Omega_+$, and with $\widetilde q$ and $\widetilde f$ defined by,
\begin{equation} \label{dataplusvsdataminus}
\widetilde q = \left\{ \begin{array}{lc} q & \mathrm{on} \ \Gamma_N \cap \partial \Omega_A \\[5pt] \dfrac{\kappa_A}{\kappa_B} q & \mathrm{in} \ \Gamma_N \cap \partial \Omega_B, \end{array} \right. \qquad \widetilde f = \left\{ 
\begin{array}{lc} f & \mathrm{in} \ \Omega_A \\[5pt] \dfrac{\kappa_A}{\kappa_B} f & \mathrm{in} \ \Omega_B. \end{array} 
\right.
\end{equation}
	 		 
   \item[$\rightarrow$]  \textsc{Local}
	\begin{eqnarray}\label{localPbm}	
&	 \Int{\Omega_-}{}{ \kappa_- \nabla T_- \cdot \nabla v}{\Omega} = \Int{\Omega_-}{}{fv}{\Omega} + 
\Int{\partial \Omega_- \cap \Gamma_N}{}{qv}{\Gamma}, \qquad \forall v \in H^1_{\gamma} (\Omega_-), \\[5pt]
& T_-(\partial \Omega_- \cap \Gamma_D) = T_D, \qquad  T_-(\gamma_i) = T_+ (\gamma_i)	
\end{eqnarray}
where  $\nn$ denotes the unit normal pointing outwards from $\Omega_-$,  and $\gamma = \gamma_i \cup (\Gamma_D \cap\partial \Omega_-)$  (for uniformity of notation we set $\kappa_- \equiv \kappa_B$).				
	\end{itemize}	
	\end{problem}

\textcolor{black}{A more detailed explanation describing the derivation of Problem 2, and establishing its consistency with Problem 1, may be found in \cite{VBA2019} }.	The following theorem establishes the consistency of the coupled formulation (\ref{globalPbm})-(\ref{localPbm}) with Problem \ref{origProblem}; it was proved in \cite{VBA2019}. Such consistency was also shown for non-constant coefficient and unsteady variants of the basic problem; however for the purposes of this work we will restrict our attention to the steady case with constant coefficients.

\noindent \textbf{Remark. } We note that $\kappa_+$ and $\kappa_-$ are quantities defined after discretization, and for this reason we elect to keep their definitions distinct from $\kappa_A$ and $\kappa_B$. In particular, $\kappa_+$ is defined on all of $\Omega_+$, while $\kappa_A$ is defined only in $\Omega_A$. For this reason, $\kappa_+$ naturally takes the defintion of $\kappa_A$ on $\Omega_A$, however, as its domain also includes $\Omega_B$, it requires an extension over this region. The implications of this extension are important from the numerical point of view and are explored in section 5.5.

\begin{theorem}\label{thm:cons} Let $\kappa_A = \kappa_+$, $\kappa_B=\kappa_-$, with $\kappa_+$ and $\kappa_-$ positive constants,  and let $\widetilde q$ and $\widetilde f$ be given by \eqref{dataplusvsdataminus}.
	If a function $T \in H^1(\Omega_+)$, \silvia{with $\Delta T|_{\Omega_-} \in L^2(\Omega_-)$ and $\Delta T|_{\Omega_+ \setminus \Omega_-} \in L^2(\Omega_+\setminus \Omega_-)$,} is a weak solution of \textcolor{black}{Problem \ref{origProblem}},  then, for $T_+ = T$, $T_- = T|_{\Omega_-}$, the couple $(T_+,T_-)$ solves  %\textcolor{black}{Problems (\ref{pb1NoAlg}).
	\textcolor{black}{Problem \ref{CoupledProblem}}.
	Conversely, if $(T_+,T_-) \in H^1(\Omega_+) \times H^1(\Omega_-)$ \silvia{with $\Delta T_+|_{\Omega_-} \in L^2(\Omega_-)$ and $\Delta T_+|_{\Omega_+ \setminus \Omega_-} \in L^2(\Omega_+\setminus \Omega_-)$,} solves  
	\textcolor{black}{Problem \ref{CoupledProblem}},	 then $T_- = T_+$ in $\Omega_-$, and $T_+$ solves  \textcolor{black}{Problem \ref{origProblem}}.

\end{theorem}

\silvia{
\begin{remark} Theorem \ref{thm:cons} was proven in \cite{VBA2019} under the stronger assumption that $T|_{\Omega_-} \in H^2(\Omega_-)$ and $T|_{\Omega_+\setminus \Omega_-} \in H^2(\Omega_+\setminus\Omega_-)$. It is however not difficult not see that, under the assumption $\bar \gamma_i \cap \Gamma_N = \emptyset$, thanks to the identity \eqref{intbyparts}, such strong assumption can be replaced by the weaker assumptions considered above.  
\end{remark}
}

\	

\noindent {\bf Two-level DD Method}. The split problem formulation (\ref{globalPbm})-(\ref{localPbm}) forms the backbone of our algorithm for solving (\ref{origProblem}):

\begin{description}	
	\item[Step 0:] Solve	
		\begin{equation}\label{eq:k0}
	\begin{array}{l}
	\kappa_+ \Int{\Omega_+}{}{ \nabla T_+^0 \cdot \nabla v}{\Omega}
	= \Int{\Omega_+}{}{fv}{\Omega} + \Int{\Gamma_N}{}{\widetilde q v}{\Gamma}, \qquad \forall v \in H^1_{\Gamma_D}(\Omega_+)\\[5pt]
	T_0(\Gamma_D) = T_D
	\end{array}
	\end{equation}
	
	\item[Step ($k+{1/2}$):] Given $T_+^k \in H^1(\Omega_+)$, compute $T_-^{k+1/2}$ by solving
	\begin{equation}\label{eq:k12}
\begin{array}{l}
\kappa_- \Int{\Omega_-}{}{ \nabla T_-^{k+1/2} \cdot \nabla v}{\Omega} = 
\Int{\Omega_-}{}{fv}{\Omega} + 
\Int{\partial \Omega_- \cap \Gamma_N}{}{qv}{\Gamma},\qquad \forall v \in
 {H^1_\gamma}(\Omega_-) \\[5pt]
\quad T_-^{k+1/2}(\partial \Omega_- \cap \Gamma_D) = T_D, \qquad 
T_-^{k+1/2}(\gamma_i) = T_+^k (\gamma_i)\end{array}
\end{equation}
where $\GmD = \Gamma_D \cap \partial\Omega_-$ and $H^1_{\gamma} = \{ v \in H^1(\Omega_-): \ v = 0 \ \text{ on } \gamma \}$.

\

	\item[Step ($k+1$):] given $T_-^{k+1/2} \in H^1(\Omega_-)$, compute $\widetilde T_+^{k+1} \in H^1(\Omega_+)$ by solving
	\begin{eqnarray}\label{eq:k1}
&\kappa_+ \Int{\Omega_+}{}{\nabla \widetilde{T}_+^{k+1} \cdot \nabla v}{\Omega} = \Int{\Omega_+}{}{fv}{\Omega} 
+ (\kappa_+ - \kappa_-)  \duality{\gamma_i}{ 
\nabla T_-^{k+1/2} \cdot \nn}{ v }
+ \Int{\Gamma_N}{}{\widetilde q v}{\Gamma}, \qquad \forall v \in H^1_{\Gamma_D}(\Omega_+),\\[5pt]
& \widetilde{T}_+^{k+1}(\Gamma_D) = T_D.
\end{eqnarray}

	\item[Relaxation step:] Set
	\begin{equation}\label{eq:k2}
T_+^{k+1} = \theta \widetilde{T}_+^{k+1} + \left(1-\theta\right) T_+^k,\qquad 0< \theta \leq 1.
\end{equation}

	\item[Check convergence:]  If convergence criteria are met, terminate, otherwise, repeat steps $k+1/2$ and $k+1$.
		
\end{description}	
	
Notice that, in Step $k+ 1/2$, the Dirichlet boundary condition $T_-^{k+1/2}(\gamma_i) = T_+^k (\gamma_i)$ can be enforced weakly (by penalization or more sophisticated approaches). Here we stick to a traditional essential treatment of the Dirichlet conditions (i.e. in the functional space). Note that  later, when we discuss the discrete problem, we will use a penalization formulation to enforce the Dirichlet-type interface condition on $\gamma_i$.
	
\noindent \textbf{Remark. } \textcolor{black}{Before continuing with the analysis, we would like to briefly discuss how the above method differs from extant methods in the literature. The Two-level DD method can be seen as a combination of two ingredients: a non-overlapping domain decomposition method (similar to \cite{da2012overlapping, toselli2004domain, R1992}) combined with a fictitious domain method \cite{ABGLR:2015, RAB2007}. For the domain decomposition, we decompose the domain into the local problem \ref{localPbm}, in which we consider a fitted mesh for $\Omega_B$. For the global problem \ref{globalPbm}, we then consider the problem on $\Omega_B$ in an unfitted manner, applying a fictitious domain discretization. The transmission conditions considered here are Dirichlet-Neumann (for the local and global problems, respectively); however, in principle, we may use any other number of transmission conditions, including Robin-Robin \cite{lions1988schwarz, lions1990schwarz, toselli2004domain} or the optimization-based approach of \cite{gunzburger1999optimization, Gunzburgeretal2000}. While the related fictitious domain and non-overlapping domain decomposition methods are well-studied in the literature, an approach combining these two paradigms, as is done here, is, to the authors knowledge, novel, though it has some strong similarities with the approach of \cite{ABGLR:2015}, the main difference being that in this last paper the transmission between local and global problem is obtained via a distributed multiplier involving possibly cumbersome integrals of the product of functions "living" on two different meshes.}

 \section{Convergence Theorem}\label{sec:conv}
 In this section we establish the convergence of the iterative procedure (\ref{eq:k0})-(\ref{eq:k2}). 
 This Theorem is one of the novel contributions of this paper.
\newline\noindent \textcolor{black}{\textbf{Theorem: } \textit{The iterative procedure defined by (\ref{eq:k0})-(\ref{eq:k2}) converges to the solution $T$ of Problem 1, provided that $\theta$ is sufficiently small for $\kappa_+/\kappa_->2$.}}
\par  We begin by analyzing the homogenous case in which $f,\,\widetilde{f},\,q,\,\widetilde{q}$, and $T_D$ are uniformly zero. 
 We let  $\boldsymbol{G}: H^1(\Omega_-) \to H^1_{\Gamma_D}(\Omega_+)$ and $ \boldsymbol{L}: H_{\Gamma_D}^1(\Omega_{+}) \to H^1_{\GmD}(\Omega_-)$ respectively denote the solution operators for the global problem (\ref{eq:k1}) and for the local problem (\ref{eq:k12}) with homogeneous data. Given $T_+^k\, \in H_{\Gamma_D}^1(\Omega_+)$, $T_+^{k+1}$ is given by
 \begin{align}
\begin{split}\label{iterator}
	T_+^{k+1} &= \theta \boldsymbol{G}\circ \boldsymbol{L} \,\lbrack T_+^k \rbrack + (1-\theta) T_+^k.
\end{split}\end{align}

We start by proving a stability bound. 
\silvia{
We observe that $v = T_-^{k+1/2} - T^k_+ \in H^1_\gamma(\Omega_-)$. We can then take such a function as a test function in (\ref{eq:k12}) and  (for $f=0, q=0$) we get
\begin{align}\begin{split}\label{localEstimate1}
	\kappa_- \int_{\Omega_-} \big|\nabla \tminus \big|^2 - \kappa_- \int_{\Omega_-} \nabla \tminus \cdot \nabla  T_+^k &= 0.
\end{split}\end{align}
}

%By letting $v = T_-^{k+1/2}$ in (\ref{eq:k12}) (for $f=0, q=0$), we observe that:
%\begin{align}\begin{split}\label{localEstimate1}
%	\kappa_- \int_{\Omega_-} \big|\nabla \tminus \big|^2 &= \bigg\lbrack  - \kappa_-\int_{\Omega_-} \left(\Delta \tminus\right)\, \tminus +  \kappa_- \int_{\partial \Omega_-} \left( \nabla \tminus\cdot \boldsymbol{n}_-\right) \tminus \bigg\rbrack.
%\end{split}\end{align}
%\noindent Since $f=0$ and $q=0$, we have that $-\Delta \tminus = 0 $ in $\Omega_-$ and $\nabla \tminus \cdot \nn=0$ on $\Gamma_N \cap \partial \Omega_- $. 
%This causes the first term  on the right hand side to vanish, and we can write:
We then immediately get the stability bound
\begin{align}\begin{split}\label{localEstimate2}
\kappa_- \big| \tminus \big|_{H^1(\Omega_-)}^2 \leq  \kappa_- \big| \tminus \big|_{H^1(\Omega_-)} \big| T_+ \big|_{H^1(\Omega_-)},
\end{split}\end{align}
from which it follows that:
\begin{align}\label{localHomogBound}
\big| \tminus \big|_{H^1(\Omega_-)} &\leq \big| T_+^k \big|_{H^1(\Omega_-)}.
\end{align}
We may then follow a similar procedure for (\ref{eq:k1}) in order to bound $\widetilde T_+^{k+1} = \boldsymbol{G} \circ \boldsymbol{L} [T_+^k]$: using standard trace bounds for harmonic functions we can write
\begin{align}\begin{split}\label{globalEstimate1}
\kappa_+ \big| \tplusw \big|_{H^1(\Omega_+)}^2 &= \left(\kappa_+-\kappa_-\right) \duality{\gamma_i} {\nabla \tminus \cdot \boldsymbol{n}} {\tplusw} \\
&\leq \big|\kappa_+-\kappa_-\big| \big|\tminus\big|_{H^1(\Omega_-)} \big| \tplusw \big|_{H^1(\Omega_-)},
\end{split}\end{align}
implying (recalling that $\kappa_+$ and $k_-$ are both positive):
\begin{align}\begin{split}\label{globalHomogBound}
	\big| \tplusw \big|_{H^1(\Omega_+)} &\leq \frac{\big| \kappa_+-\kappa_-\big|}{%\big| 
	\kappa_+ %\big|
} \big| \tminus \big|_{H^1(\Omega_-)} \\
	&\leq \frac{\big| \kappa_+-\kappa_-\big|}{%\big| 
	\kappa_+ %\big|
}  \big| T_+^k \big|_{H^1(\Omega_+)}.
\end{split}\end{align}
Observe, at this point, that, if
\begin{equation}\label{cond30}
\frac{| \kappa_+ - \kappa_- |}{ \kappa_+ } < 1 \Leftrightarrow (0 < ) \ \dfrac{\kappa_-}{2} < \kappa_+ 
\end{equation}
then the operator $\boldsymbol{G}\circ \boldsymbol{L}$ is a contraction, and the sequence $T_+^k$ converges to a unique fixed point for all $\theta \in ]0,1]$. In particular, we can then take $\theta=1$ and the relaxation step (\ref{eq:k2}) is not necessary.

%If $\kappa_+ \geq \kappa_-$, then we easily check that \eqref{cond30} always holds. 
We then consider the case $\kappa_+ \leq \kappa_-/2$. We let $\tplusw = \boldsymbol{G} \circ \boldsymbol{L}\, \lbrack T_+^k \rbrack$ and recall that $f,\,\widetilde{f},\,q,\,\widetilde{q}$ and $T_D$ are all zero. We can write
	\begin{align}\begin{split}\label{case2Bound1}
\kappa_+ \int_{\Omega_+} \nabla \tplusw \cdot \nabla T_+^k &=  (\kappa_+-\kappa_-)\duality{\gamma_i}{\nabla \tminus \cdot \nn} {T_+^k} \\ &= (\kappa_+ - \kappa_-) \duality{\gamma_i} {\nabla \tminus \cdot \nn}{ \tminus} \\
&= (\kappa_+ - \kappa_-)  \duality{\partial \Omega_-} {\nabla \tminus \cdot \nn} {\tminus},\end{split}\end{align}
where we exploit the fact that $T^{k+1/2}_- = T^k_+$ on $\gamma_i$ and
the homogeneous boundary conditions on $\Gamma_D,\, \Gamma_N$, where
$\partial \Omega_- = \gamma_i \cup (\partial \Omega_- \cap \partial \Omega_+)$.  Integration by parts, $ \Delta \tminus = 0$ in $\Omega_-$, and $(\kappa_+ - \kappa_- )< 0$ yield:
  	\begin{align}\begin{split}\label{case2Bound2}
\kappa_+ \int_{\Omega_+} \nabla \tplusw \cdot \nabla T_+^k &= (\kappa_+-\kappa_-) \left\lbrack  \int_{ \Omega_-} ( \Delta \tminus )\tminus + \int_{\Omega_-}  \big| \nabla \tminus \big|^2 \right\rbrack \\
&= (\kappa_+-\kappa_-) | \tminus |_{H^1(\Omega_-)}^2 \leq 0.
\end{split}\end{align}
\par We then apply relaxation as in (\ref{eq:k2}), giving:
\begin{equation}\label{case2Relax}
\tplus = \theta \tplusw + (1-\theta)T_+^k = T_+^k + \theta(\tplusw - T_+^k).	
\end{equation} 
It then follows that:
\begin{align}\begin{split}\label{case2FinalIneq}
\kappa_+ \big| \tplus \big|_{H^1(\Omega_+)}^2 &= \kappa_+\bigg\lbrack \big| T_+^k \big|_{H^1(\Omega_+)}^2 + \theta^2 \big| \tplusw \big|_{H^1(\Omega_+)}^2 + \theta^2 \big| T_+^k \big|_{H^1(\Omega_+)}^2  \\
 &\quad - 2\theta^2(\nabla \tplusw,\,\nabla T_+^k) + 2\theta(\nabla \tplusw,\,\nabla T_+^k) -2\theta \big| T_+^k\big|_{H^1(\Omega_+)}^2\bigg\rbrack \\
 &= \kappa_+\left(1+\theta^2-2\theta\right)\big| T_+^k \big|_{H^1(\Omega_+)}^2 + \kappa_+ \theta^2 \big| \tplusw \big|_{H^1(\Omega_+)}^2 \\
 &\quad + \kappa_+ 2(\theta-\theta^2) (\nabla \tplusw, \nabla T_+^k) \\
 &\leq  \kappa_+\left(1-\theta\right)^2\big| T_+^k \big|_{H^1(\Omega_+)}^2  +  
 \theta^2 \dfrac{|\kappa_+-\kappa_-|^2}{\kappa_+} \ \big| T_+^k \big|_{H^1(\Omega_+)}^2,
\end{split}\end{align}
where the last line follows from (\ref{globalHomogBound}) and (\ref{case2Bound2}). 
From (\ref{case2FinalIneq}), and the fact that $| \kappa_+ - \kappa_-|$ in this case is $\kappa_- - \kappa_+$:
\begin{align}\label{finalIneqHomog}
\big| \tplus \big|_{H^1(\Omega_+)}^2 &
\leq \left(
 (1 - \theta)^2 + \theta^2 \left( \dfrac{(\kappa_- -\kappa_+)^2}{\kappa_+^2} \right) \right) \big|T_+^k\big|_{H^1(\Omega_+)}^2, 
%\theta^2\left(1+\frac{|\kappa_+ - \kappa_-|}{\kappa_+}\right) - 2\theta \right) \big|T_+^k\big|_{H^1(\Omega_+)}^2
\end{align}
so that
\begin{align}\label{finalIneqHomog2}
	\big| \tplus \big|_{H^1(\Omega_+)}^2 &\leq \left( \left(\left(\dfrac{\kappa_- - \kappa_+}{\kappa_+} \right)^2 +1 \right) \theta^2 - 2 \theta + 1 \right) 
\big|T_+^k\big|_{H^1(\Omega_+)}^2
	%\theta^2\left(1+\frac{|\kappa_+ - \kappa_-|}{\kappa_+}\right) - 2\theta \right) \big|T_+^k\big|_{H^1(\Omega_+)}^2
\end{align}

The parabola $ \left(\left(\dfrac{\kappa_- - \kappa_+}{\kappa_+}\right)^2 +1 \right) \theta^2 - 2 \theta + 1$ has value 1
and slope negative for $\theta=0$ (the slope being -2), minimum in $\theta_{opt} = \left(\left(\dfrac{\kappa_- - \kappa_+}{\kappa_+}\right)^2 +1 \right)^{-1}$ 
with value $1 - \theta_{opt} < 1$. So, there exists an interval $(0,\bar{\theta}$) where it takes values
$< 1$, proving that the map is a contraction also in this case.

%One then verifies that for $\theta < \frac{2 \kappa_+}{\kappa_-}$, 
%\begin{align}\label{globalHomogContractCase2}
%\left(1+ \theta^2\left(1+\frac{|\kappa_+-\kappa_-|}{\kappa_+}\right) - 2\theta \right) <1, \end{align}
%so that, if $k_+ < k_-$, provided the relaxation parameter $\theta$ verifies $\theta < 2k_+/k_-$, the operator mapping $T_+^k$  to $T_+^{k+1}$ is a contraction, and the iterative procedure converges to a unique fixed point. Remark that if $k_+< \kappa_-$ verifies $\kappa_+ > k_-/2$, then  $\theta = 1$ verifies $\theta < 2\kappa_+ /\kappa_-$ and, also in such case, no relaxation is needed.

\
 
 \par We now consider the non-homogeneous case, with $f,\,\widetilde{f} \in L^2(\Omega_+) $ and $q,\,\widetilde{q} \in H^{-1/2} (\Gamma_N) $ and $T_D \in H^{1/2}(\Gamma_D)$. Let $\tlaststar,\,\tpluswstar \in H^1(\Omega_+)$ and $\tminusstar \in H^1(\Omega_-)$ be auxiliary functions defined as the solutions to the following problems, the first one being a standard elliptic boundary value problem with given data, the second and third being also standard elliptic boundary value problems with data depending on the solution of the first and second problem, respectively:
 \begin{align}\begin{split}\label{TKstar}
	\int_{\Omega_+} \kappa_+ \nabla \tlaststar \cdot \nabla v &= \int_{\Omega_+} \widetilde{f} v + \int_{\Gamma_N} \widetilde{q} v,  \quad \forall v \in H_{\Gamma_D}^1(\Omega_+), \qquad
	 \tlaststar(\Gamma_D)=T_D;
	\end{split} 
		\end{align}

 \begin{align}
 \begin{split}
 \label{TMinusStar}
	&\int_{\Omega_-} \kappa_- \nabla \tminusstar \cdot \nabla v =  \int_{\Omega_-} f v + \int_{\Gamma_N \cap \partial \Omega_- } q v,\qquad  \forall v \in H_{\gamma }^1(\Omega_-) \\ &\tminusstar(\Gamma_D \cap \partial \Omega_- )=T_D, \qquad \tminusstar(\gamma_i)=\tlaststar)(\gamma_i) ;
	\end{split} \end{align}
	
	 \begin{align}\begin{split}\label{TPlusstar}
	\int_{\Omega_+} \kappa_+ \nabla \tpluswstar \cdot \nabla v &=  (\kappa_+-\kappa_-)\duality{\gamma_i}  {\nabla \tminusstar \cdot \nn}{v},\qquad   \forall v \in H_{\Gamma_D}^1(\Omega_+), \qquad \tpluswstar(\Gamma_D)=0.
	\end{split} \end{align}
	
	To study the convergence of the sequence $T_+^k$ of the iteration \eqref{eq:k0}--\eqref{eq:k2}, we investigate the existence of a limit to the sequence
	 \[ w_+^k = T_+^k - T_+^*.\] To this aim we show that the operator mapping $w_+^k$ to $w_+^{k+1}$ is a contraction.
	We start by rewriting (\ref{eq:k12}) and (\ref{eq:k1}) in terms of $w_+^k$ and of the new local unknown \[w_-^{k+1/2} = T_-^{k+1/2} - T_-^*.\]
	This gives us the following equation for $w_-^{k+1/2}$: 
 	\begin{equation}\label{eq:k1Hom}
\begin{array}{l}
\kappa_- \Int{\Omega_-}{}{ \nabla (w_-^{k+1/2} + \tminusstar) \cdot \nabla v}{\Omega}   = 
\Int{\Omega_-}{}{fv}{\Omega} + 
\Int{\partial \Omega_- \cap \Gamma_N}{}{qv}{\Gamma}, \qquad \forall v \in H^1_{\Gamma_D \cup \gamma_i} (\Omega_-)\\[5pt]
\quad (\wminus + \tminusstar )(\partial \Omega_- \cap \Gamma_D) = T_D, \quad 
(\wminus + \tminusstar )(\gamma_i) = (\wlast + \tlaststar)(\gamma_i).
\end{array}
\end{equation}
Substituting (\ref{TMinusStar}) into (\ref{eq:k1Hom}), it follows promptly that:
\begin{align}\label{nonHomogIntoHomogLocal}
	\wminus &= \boldsymbol{L} \wlast,
\end{align}
where $\boldsymbol{L}$ is the homogeneous local-problem solution operator studied above.
\par Following a similar procedure for (\ref{eq:k1}), we see that, setting $\wplus = \widetilde T_+^{k+1} - T_+^*$, we can write, for all $v$ in $H^1_{\Gamma_D}(\Omega_+)$,
\begin{align}\begin{split}
	\int_{\Omega_+} \kappa_+ \nabla\left(\wplus + \tplusstar \right)\cdot\nabla v &= \int_{\Omega_+} \widetilde{f}v + \int_{\Gamma_N} \widetilde{q}v+ (\kappa_+ - \kappa_-) \duality{\gamma_i}{\nabla \tminus\cdot \nn}{v}, \qquad \forall v \in H^1_{\Gamma_D}(\Omega_+)\\
	& \wplus(\Gamma_D) + T_+^*(\Gamma_D) = T_D
\end{split}\end{align}
whence, using \eqref{TKstar} and \eqref{TPlusstar}
\begin{align}\begin{split}\label{nonHomogIntoGlobal}
		\int_{\Omega_+} \kappa_+ \nabla \wplus \cdot\nabla v &= (\kappa_+ - \kappa_-) \duality{\gamma_i}{\nabla \left( \wminus+\tminusstar\right) \cdot \nn} {v} \\
			&= (\kappa_+ - \kappa_-) \duality{\gamma_i}{\nabla \left( \wminus \right) \cdot \nn} {v}  + \int_{\Omega_+}\kappa_+ \nabla \tpluswstar \cdot \nabla v.
\end{split}\end{align} 
Then, $\wplus - \tpluswstar = \boldsymbol{G} \wminus$, and, therefore, 
\begin{equation}\label{wIterator}
	\wplus = \boldsymbol{G}\lbrack \wminus \rbrack + \tpluswstar = \boldsymbol{G} \circ \boldsymbol{L} \lbrack w_+^k \rbrack + \tpluswstar,
\end{equation}
 Applying relaxation, one obtains:
\begin{align}\begin{split}\label{conclusion}
	w_+^{k+1} &= \theta \widetilde{w}_+^{k+1} + (1-\theta)w_+^k \\
	&= \theta(\boldsymbol{G} \circ \boldsymbol{L} \lbrack w_+^k \rbrack  + \tpluswstar ) + (1-\theta)w_+^k \\
	& = \left(\theta \boldsymbol{G} \circ \boldsymbol{L} \lbrack w_+^k \rbrack  + (1-\theta)w_+^k \right) + \theta \tpluswstar.
\end{split}\end{align}
%
%
%ALEVEN: HO CAMBIATO QUI
%
%We refer to \eqref{conclusion} hereafter, assuming that $\theta = 1$ for
%$\kappa_-/\kappa_+ < 2$,  while in the other case we assume $\theta$ to be in the range 
%of contractivity.
%
From \eqref{conclusion}, the map $w_+^{k+1} = {\cal M} \lbrack w_+^{k+1} \rbrack$
is contractive. 

It follows that $\tplus = \tlaststar + w_+^{k+1}$, obtained directly by solving the nonhomogeneous problem, converges to the solution of the coupled problem (\ref{globalPbm})-(\ref{localPbm}).

 \begin{remark}
 	Observe that, in order for the iterative procedure to converge, we did not require that $q$ and $\widetilde q$ (resp. $f$ and $\widetilde f$) satisfy \eqref{dataplusvsdataminus}. Of course, such relations are needed if we want the solution of the coupled problem to coincide with the solution of the original problem, as stated in Theorem \ref{thm:cons}. \silvia{
 	However, it could be interesting to exploit the freedom in the choice of the extension of $f$ in the design of the method. In particular, in the spirit of \cite{Atamianetal:1991}, one could look for an extension $\widetilde f$ such that the jump along $\gamma_i$ of $\nabla T^+ \cdot \nn$ vanishes, thus allowing for optimal convergence rates.
 	}
 	\end{remark}

\section{Discrete Problem}\label{sec:discrete}
Having formally established the convergence of the continuous problem in the preceding section, in the present section we consider the discrete version of Problem (\ref{origProblem}) and of its corresponding Two-level DD formulation (\ref{eq:k0})-(\ref{eq:k1}). We first introduce a discrete \textit{monolithic problem}. Let $\mathcal{T}_+$ and $\mathcal{T}_-$ denote appropriate discretizations of $\Omega_+$ and $\Omega_-$, \silvia{with $\Edgesi$ denoting the discretization of $\gamma_i$ induced by $\mathcal{T}_-$,} and let $X_{\mathcal{T}_+} \subset H^1(\Omega_+)$, $X_{\mathcal{T}_-} \subset H^1(\Omega_-)$ denote the corresponding finite dimensional approximation spaces. Using a penalization method with parameter $\alpha > 0$  for the local problem to impose the Dirichlet-type coupling condition on $\gamma_i$, we consider the following discrete equations:
 find $(T_+^h, T_-^h)$ in $X_{\mathcal{T}_+}  \times X_{\mathcal{T}_-}$ such that for all $v_h$ in $X_{\mathcal{T}_+} \cap H^1_{\Gamma_D}(\Omega_+)$ and $w_h$ in $X_{\mathcal{T}_-} \cap H^1_{\Gamma_{-\,D}}(\Omega_-)$:
\begin{align}\begin{split}\label{monolithicPbm}
\int_{\Omega_+} \kappa_+ \nabla T_+^h \cdot \nabla v_h - \silvia{\sum_{e \in \Edgesi} \int_{e} \left(\kappa_+ -\kappa_-\right) \left(\nabla T_-^h \cdot \nn\right)v_h} &= \int_{\Omega_+} \widetilde{f}\,v_h + \int_{\Gamma_N\cap \partial \Omega_-}\widetilde{q}\,v_h \\
\int_{\Omega_-} \kappa_- \nabla T_-^h \cdot \nabla w_h + \alpha  \int_{\gamma_i}T_-^h w_h - \alpha  \int_{\gamma_i}T_+^h w_h &= \int_{\Omega_-} f\,w_h + \int_{\Gamma_N \cap \partial \Omega_- } q\,w_h \\
T_+^h &= T_0 \text{ on } \Gamma_D.
	\end{split}\end{align}
 respectively, with corresponding discrete function spaces $\XTO$ and $\XTT$. We denote discrete functions with the subscript $h$. 
\par  We define the following matrices resulting from the discretization of the bilinear forms  in (\ref{monolithicPbm}):
\begin{align}\label{KGlobal}
	K_+ & \ \text{for}  \ \int_{\Omega_+} \kappa_+ \nabla T_+^h \cdot \nabla v_h 
\end{align}
\begin{align}\label{gammaGlobal}
	S & \ \text{for}  \  \silvia{\sum_{e \in \Edgesi} \int_{e} \left(\kappa_+ -\kappa_-\right) \left(\nabla T_-^h \cdot \nn\right)v_h}
\end{align}
\begin{align}\label{diricLocal}
	D &\ \text{for}  \  -\alpha \int_{\gamma_i} T_+^h w_h
\end{align}
\begin{align}\label{KLocal}
	K_- & \ \text{for}  \  \int_{\Omega_-} \kappa_- \nabla T_-^h \cdot \nabla w_h + \alpha \int_{\gamma_i} T_-^h w_h,
\end{align}
using which the algebraic form of (\ref{monolithicPbm}) reads
\begin{align}\label{algSyst}
\begin{bmatrix}
K_+ & S  \\ D & K_-	
\end{bmatrix}
	\begin{bmatrix}
	\boldsymbol{T}_+^h \\ \boldsymbol{T}_-^h	
	\end{bmatrix}
&=
\begin{bmatrix}
\boldsymbol{f}_+ \\ \boldsymbol{f}_-	
\end{bmatrix}.\end{align}
We now offer two useful and equivalent algebraic interpretations of Algorithm in (\ref{eq:k0})-(\ref{eq:k2}).

\subsection{Gauss-Seidel-type interpretation}
We note that (\ref{monolithicPbm}) can be written equivalently as: 
\begin{align}\label{additiveSplit}
\begin{bmatrix}
K_+ & S  \\ 0 & K_-	
\end{bmatrix}
	\begin{bmatrix}
	\boldsymbol{T}_+^h \\ \boldsymbol{T}_-^h	
	\end{bmatrix}
&=\begin{bmatrix}
\boldsymbol{f}_+ \\ \boldsymbol{f}_-	
\end{bmatrix} 
+\begin{bmatrix}
0 & 0  \\ -D & 0	
\end{bmatrix}
	\begin{bmatrix}
	\boldsymbol{T}_+^h \\ \boldsymbol{T}_-^h	
	\end{bmatrix}. \end{align}
The above splitting (\ref{additiveSplit}) can then be employed to solve (\ref{algSyst}) iteratively as an inverted \textit{block Gauss-Seidel}-type method\footnote{The reason why we call it ``inverted'' is because it works with an upper triangular matrix as opposed to a lower triangular one, as in the standard GS.}- given $\lbrack \boldsymbol{T}_+^{h,0},\, \boldsymbol{T}_-^{h,0}\rbrack^T$, one solves until convergence:
\begin{align}\label{gaussSeidel}
\begin{bmatrix}
K_+ & S  \\ 0 & K_-	
\end{bmatrix}
	\begin{bmatrix}
	\boldsymbol{T}_+^{h,k+1} \\ \boldsymbol{T}_-^{h,k+1}	
	\end{bmatrix}
&=\begin{bmatrix}
\boldsymbol{f}_+ \\ \boldsymbol{f}_-	
\end{bmatrix} 
+\begin{bmatrix}
0 & 0  \\ -D & 0	
\end{bmatrix}
	\begin{bmatrix}
	\boldsymbol{T}_+^{h,k} \\ \boldsymbol{T}_-^{h,k}	
	\end{bmatrix}. \end{align}
One may modify (\ref{gaussSeidel}) to incorporate an under-relaxation parameter $\theta$,
\begin{align}\label{gaussSeidelRelax}
\begin{bmatrix}
\frac{1}{\theta} K_+ & S  \\ 0 & K_-	
\end{bmatrix}
	\begin{bmatrix}
	\boldsymbol{T}_+^{h,k+1} \\ \boldsymbol{T}_-^{h,k+1}	
	\end{bmatrix}
&=\begin{bmatrix}
\boldsymbol{f}_+ \\ \boldsymbol{f}_-	
\end{bmatrix} 
+\begin{bmatrix}
\frac{1-\theta}{\theta}K_+ & 0  \\ -D & 0	
\end{bmatrix}
	\begin{bmatrix}
	\boldsymbol{T}_+^{h,k} \\ \boldsymbol{T}_-^{h,k}	
	\end{bmatrix}. \end{align} 
The scheme (\ref{gaussSeidel}) leads to solving the sequence of problems:
\begin{align}\begin{split}\label{algLocalPbm}
K_- \boldsymbol{T}_-^{h,k+1} &= \boldsymbol{f}_- - D \boldsymbol{T}_+^{h,k}	\,,
\end{split}\\
\begin{split}\label{algGlobalPbm}
K_+ \boldsymbol{T}_+^{h,k+1} &= \boldsymbol{f}_+ - S \boldsymbol{T}_-^{h,k+1}\, ,
\end{split}\end{align}
where 
(\ref{algLocalPbm}) and (\ref{algGlobalPbm}) are the algebraic versions of step $k+1/2$ and $k+1$, respectively, in (\ref{eq:k0})-(\ref{eq:k1}). Incorporating the  under-relaxation step (\ref{eq:k2}) is equivalent to solving the modified system (\ref{gaussSeidelRelax}):
\begin{align}\begin{split}\label{algLocalPbmRelax}
K_- \boldsymbol{T}_-^{h,k+1} &= \boldsymbol{f}_-  - D \boldsymbol{T}_+^{h,k}	\,,
\end{split}\\
\begin{split}\label{algGlobalPbmRelax}
\frac{1}{\theta}K_+ \boldsymbol{T}_+^{h,k+1} &= \boldsymbol{f}_+  + \frac{1-\theta}{\theta}K_+\boldsymbol{T}_+^{h,k} - S \boldsymbol{T}_-^{h,k+1}\,.
\end{split}\end{align}
It is well-known that the block Gauss-Seidel iteration (\ref{gaussSeidelRelax}) converges to the solution of (\ref{algSyst}) provided that the spectral radius of the system:
\begin{align}\label{specRad}
\begin{bmatrix}
\frac{1}{\theta} K_+ & S  \\ 0 & K_-	
\end{bmatrix}^{-1}
\begin{bmatrix}
\frac{1-\theta}{\theta}K_+ & 0  \\ -D & 0	
\end{bmatrix}
&= \begin{bmatrix} \left(1-\theta\right)I_+ + \theta K_+^{-1}SK_-^{-1}D  & 0 \\ -K_-^{-1}D & 0   
\end{bmatrix}
\end{align}
is less than one (see e.g. \cite{GoluVanl2013}), where $I_+$ denotes the identity matrix of appropriate size for the function space $\XTO$ ($I_-$ will be used similarly). As (\ref{specRad}) features a zero block on the diagonal, the convergence criterion
reduces to the condition:
\begin{align}\label{discConvCriterion}
\rho\left( 	\left(1-\theta\right)I_+ + \theta K_+^{-1}SK_-^{-1}D \right) < 1.
\end{align}
%If we assume for simplicity that the relaxation parameter $\theta=1$, this yet further reduces to:
%\begin{align}\label{discConvCriterionThetaOne}
%\rho\left( 	 K_+^{-1}SK_-^{-1}D \right) < 1.
%\end{align}
This is the algebraic counterpart of the condition (\ref{cond30}); both rely on the relative values of $\kappa_-$ and $\kappa_+$ (as we will confirm with the numerical experiments).
\subsection{Neumann series interpretation ($\theta = 1$)}
Following classical arguments related to the Schur complement \cite{GoluVanl2013},
we may factor (\ref{algSyst}) in the following block-UL system:
\begin{align}\label{blockUL}
\begin{bmatrix}
K_+ - SK_-^{-1}D & S  \\ 0 & K_-	
\end{bmatrix}
\begin{bmatrix}
I_{+} & 0 \\ K_-^{-1}D & I_{-} 	
\end{bmatrix}
	\begin{bmatrix}
	\boldsymbol{T}_+^h \\ \boldsymbol{T}_-^h	
	\end{bmatrix}
&=
\begin{bmatrix}
\boldsymbol{f}_+ \\ \boldsymbol{f}_-	
\end{bmatrix}.\end{align}
Solving the system yields:
\begin{align}\begin{split}\label{expandedSol}
	\Tpl &= \left\lbrack K_+ - SK_-^{-1}D\right\rbrack^{-1}\left(\fp - SK_-^{-1}\fm\right) \\
	\Tmi &= K_-^{-1} \fm - K_-^{-1}D\Tpl \\
	&= K_-^{-1} \fm - K_-^{-1}D\left\lbrack K_+ - SK_-^{-1}D\right\rbrack^{-1}\left(\fp - SK_-^{-1}\fm\right).
\end{split}\end{align}
Similarly, for an iteration $k+1$ of (\ref{algLocalPbm}) -(\ref{algGlobalPbm}) (assuming $\theta=1$):
 \begin{align}\begin{split}\label{expandedAlgTwoLevel}
	\boldsymbol{T}_+^{h,k+1} &= 	K_+^{-1} \fp - K_+^{-1}S \boldsymbol{T}_-^{h,k+1} \\
&=	K_+^{-1} \fp - K_+^{-1}S\left\lbrack K_-^{-1}\fm - K_-^{-1}D\boldsymbol{T}_+^{h,k} \right\rbrack \\
	&= K_+^{-1}\left(\fp - SK_-^{-1}\fm\right) +K_+^{-1}SK_-^{-1}D\boldsymbol{T}_+^{h,k} \\
	\boldsymbol{T}_-^{h,k+1} &= K_-^{-1}\fm - K_-^{-1}D\boldsymbol{T}_+^{h,k}.
\end{split}\end{align}
\begin{prop}\label{prop:neumann} Given an initial guess $\boldsymbol{T}_+^{h,0}$, at an iteration $k$ of (\ref{algLocalPbm})-(\ref{algGlobalPbm}) we may write $\boldsymbol{T}_+^{h,k}$ as:
$$ \boldsymbol{T}_+^{h,k} =\left\lbrack \sum_{j=0}^{k-1} \left(K_+^{-1}SK_-^{-1}D\right)^j\right\rbrack K_+^{-1} \left(\fp - SK_-^{-1}\fm\right) + \left(K_+^{-1}SK_-^{-1}D\right)^k \boldsymbol{T}_+^{h,0}.  $$ 
Further, if (\ref{discConvCriterion}) holds, then the iteration given by (\ref{algLocalPbm})-(\ref{algGlobalPbm}) converges to the solution of (\ref{algSyst}).\end{prop}
 \textbf{Proof. } Proceed by induction.
 The case $k=1$ is trivially true.
% First we verify that this holds for the base case. At $k=1$ we have that:
% \begin{align}\begin{split}\label{baseCase}
%\boldsymbol{T}_{h,1}_+  	&= K_+^{-1}\left(\fp - SK_-^{-1}\fm\right) %+K_+^{-1}SK_-^{-1}D\boldsymbol{T}_+^{h,0} 	\\
%&= \left\lbrack \sum_{j=0}^{0} \left(K_+^{-1}SK_-^{-1}D\right)^j\right\rbrack K_+^{-1} \left(\fp - %SK_-^{-1}\fm\right) + \left(K_+^{-1}SK_-^{-1}D\right)^1 \boldsymbol{T}_+^{h,0} \\
%&= \left\lbrack \sum_{j=0}^{k-1} \left(K_+^{-1}SK_-^{-1}D\right)^j\right\rbrack K_+^{-1} \left(\fp - %SK_-^{-1}\fm\right) + \left(K_+^{-1}SK_-^{-1}D\right)^k \boldsymbol{T}_+^{h,0}
%\end{split}\end{align}
Assume that the hypothesis holds for $k$. At $k+1$:
 
\begin{align}\begin{split}\label{twoLevelNeumannStep}
\boldsymbol{T}_+^{h,k+1}  	&= K_+^{-1}\left(\fp - SK_-^{-1}\fm\right) +K_+^{-1}SK_-^{-1}D\boldsymbol{T}_+^{h,k} \\
&= K_+^{-1}\left(\fp - SK_-^{-1}\fm\right) + \\
 &\quad K_+^{-1}SK_-^{-1}D\left\lbrack \left\lbrack \sum_{j=0}^{k-1} \left(K_+^{-1}SK_-^{-1}D\right)^j\right\rbrack K_+^{-1} \left(\fp - SK_-^{-1}\fm\right) + \left(K_+^{-1}SK_-^{-1}D\right)^k \boldsymbol{T}_+^{h,0} \right\rbrack \\
 &= K_+^{-1}\left(\fp - SK_-^{-1}\fm\right) + \\
 &\quad \left\lbrack \sum_{j=1}^{k} \left(K_+^{-1}SK_-^{-1}D\right)^j\right\rbrack K_+^{-1} \left(\fp - SK_-^{-1}\fm\right) + \left(K_+^{-1}SK_-^{-1}D\right)^{k+1} \boldsymbol{T}_+^{h,0} \\
 &= \left\lbrack \sum_{j=0}^{k} \left(K_+^{-1}SK_-^{-1}D\right)^j\right\rbrack K_+^{-1} \left(\fp - SK_-^{-1}\fm\right) + \left(K_+^{-1}SK_-^{-1}D\right)^{k+1} \boldsymbol{T}_+^{h,0},
 \end{split}\end{align}
completing the first part of the proof. 
\par It remains to show that if (\ref{discConvCriterion}) for $\theta=1$ holds, the solutions (\ref{expandedSol}) and (\ref{expandedAlgTwoLevel}) are equivalent as $k \to \infty$. Note that
for $\theta=1$ (\ref{discConvCriterion}) implies the convergence of the Neumann series \begin{align}\begin{split}\label{neumansSeries}
\left(K_+ - SK_-^{-1}D\right)^{-1} &= \left\lbrack K_+\left(I_+ - K_+^{-1}SK_-^{-1}D\right)\right\rbrack^{-1}	\\
&= \left(I_+ - K_+^{-1}SK_-^{-1}D\right)^{-1} K_+^{-1} \\
&=\left\lbrack \sum_{j=0}^{\infty} \left(K_+^{-1}SK_-^{-1}D\right)^j\right\rbrack K_+^{-1},
\end{split}\end{align}
hence from (\ref{expandedSol}):
\begin{align}\label{monolithicSolutionNeumann}
	\Tpl &=\left\lbrack \sum_{j=0}^{\infty} \left(K_+^{-1}SK_-^{-1}D\right)^j\right\rbrack K_+^{-1} \left(\fp - SK_-^{-1}\fm\right).
	\end{align}
Condition (\ref{discConvCriterion}) for $\theta=1$ implies additionally that:
\begin{align}\label{nilpotency}
	\lim_{k \to \infty} \left(K_+^{-1}SK_-^{-1}D\right)^k &= 0,
\end{align}
with $0$ denoting the zero matrix of appropriate dimension. Taken together, (\ref{twoLevelNeumannStep}), (\ref{monolithicSolutionNeumann}), and (\ref{nilpotency}) imply the result.

\par Proposition \ref{prop:neumann} provides powerful intuition regarding the performance of the discrete Algorithm (\ref{eq:k0})-(\ref{eq:k2}), as it can be seen as a finite approximation of the Neumann series expansion. We therefore expect that the convergence rate will depend on $\rho\left(K_+^{-1}SK_-^{-1}D\right)$. 
\subsection{Eigenvalue scaling}
Both the Gauss-Seidel (for $\theta=1$) and Neumann series interpretations of Algorithm in (\ref{eq:k0})-(\ref{eq:k2}) require the same condition for convergence: $\rho\left(K_+^{-1}SK_-^{-1}D\right) < 1$. Although a fully rigorous spectral analysis of this system is a worthy subject of future work, we will provide a more heuristic approach here which we will validate with numerical tests. 
\par Looking at the operators (\ref{KGlobal})-(\ref{KLocal}), one may argue that $\rho$ scales in the following way:
\begin{align}\begin{split}\label{scaling}
	\rho\left(K_+^{-1} S K_-^{-1} D\right) &\sim \frac{1}{\kappa_+}\left|-\kappa_+ + \kappa_-\right| \frac{1}{\kappa_- + \alpha} \alpha \\
	&= \frac{|\kappa_- -\kappa_+|}{\kappa_+}\frac{\alpha}{\kappa_-+\alpha}  \\
	&= \frac{\left|\kappa_- -\kappa_+\right|\alpha}{\kappa_+\left(\frac{\kappa_-}{\alpha} + 1\right)\alpha} \\
	&= \frac{|\kappa_- - \kappa_+|}{\kappa_+\left(\frac{\kappa_-}{\alpha}+1\right)}.
\end{split}\end{align}
We recall that $\alpha$ is a penalization parameter and is large in general, implying that $\kappa_-/\alpha \approx 0$, and therefore that:
\begin{align}\begin{split}\label{scaling2}
	\rho\left(K_+^{-1} S K_-^{-1} D\right) &\sim \frac{|\kappa_- - \kappa_+|}{\kappa_+} 
%	&= \frac{\left|\frac{\kappa_-}{\kappa_+} -1\right|\kappa_+}{\kappa_+} \\
	= \left| \frac{\kappa_-}{\kappa_+}-1 \right|.
\end{split}\end{align}

This analysis is consistent with the results of the unrelaxed scheme in Sect. \ref{sec:conv},
and specifically with \eqref{cond30}, for $\boldsymbol{G}\circ \boldsymbol{L}$ to be a contraction.

We therefore postulate the following scaling behavior:
\begin{align}\label{scalingFn}
\rho\left(\frac{\kappa_-}{\kappa_+}\right) &\sim  C_{h_+,h_-,\gamma_i ,\Omega_+,\Omega_-} \left|\frac{\kappa_-}{\kappa_+} - 1\right|,
\end{align}
where $C_{h_+,h_-,\gamma,\Omega_+,\Omega_-}$ is a constant depending on the geometric and mesh parameters and will be denoted by $\widetilde{C}$ hereafter. 

\section{Numerical experiments}\label{sec:numexp}
In this section, we will perform a series of numerical experiments to confirm the theory discussed in previous sections. We wish to confirm/examine the following:
\begin{enumerate}
\item That the convergence of the Two-level DD method is dependent on $\specRad$ as postulated;
\item That the scaling behavior (\ref{scalingFn}) holds;
\item How $\specRad$ is affected by changes in global mesh level $h_+$, local mesh level $h_-$, spatial dimension, and polynomial degree.
\end{enumerate}
To answer these questions, we will run series of 2D and 3D tests, organizing the analysis of the results analysis into two distinct portions:
\begin{enumerate}
\item \textbf{Convergence analysis.} This will focus directly on the relationship between $\specRad$ and convergence of the Two-level DD method. For a given case, we will use the collected data examine the impact of $\specRad$ on convergence. We then verify the point at which the Two-level DD method no longer converges, where $\specRad>1$.
\item \textbf{Spectral growth analysis.} After confirming the importance of $\specRad$ with the convergence analysis, we will examine the full range of results in 2D and 3D in order to understand how differences in global and local mesh resolution, degree of polynomial approximation, and spatial dimension affect the value of $\specRad$. We hope to observe that $\specRad \sim \mathcal{O}\left( \log h_+/h_-\right)$, for all polynomial degrees in both 2D and 3D, which will indicate that refinement of the local mesh does not lead to a deterioration of convergence behavior.
\end{enumerate}
The underlying data and problem setup is identical for both analyses. For this reason, we outline the setup and experimental procedure below for each case before proceeding to the discussion of results. 
\subsection{Problem setup}
\begin{figure}
\begin{center}
\includegraphics[width=.5\textwidth]{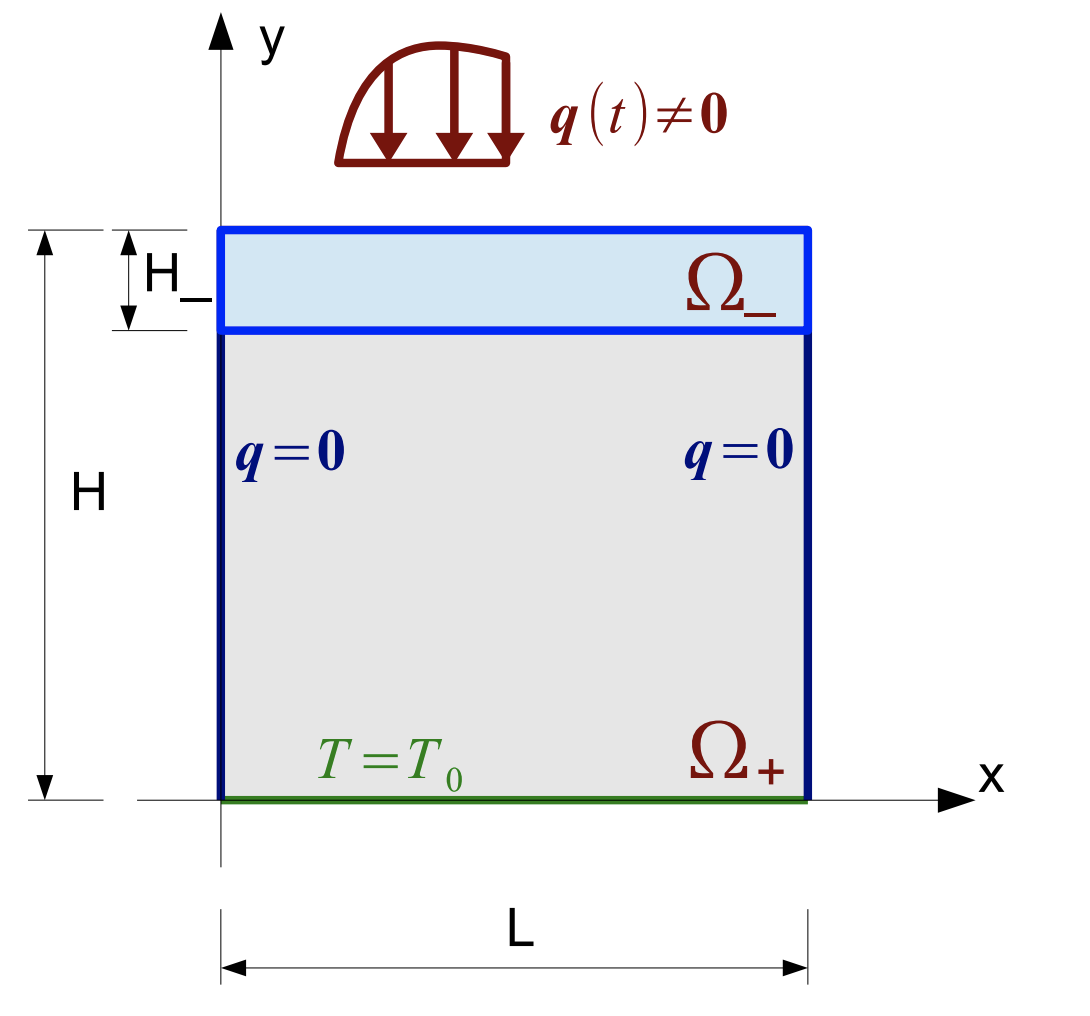}\caption{Basic test setup for the problem solved throughout this work. This is the two-dimensional version; the three-dimensional problem features the same general setup with a cubic domain rather than a square domain.}\label{fig:modelProblem}
\end{center}
\end{figure}

In order to control the experimental setting as much as possible, we define a standardized problem setup. Referring to notation shown in Fig \ref{fig:modelProblem}, we solve the 2D problem in a  square with $L=1/40$ and $H=1/40$. We define $H_-=1/160$ and set $T_0=293.15$. The heat flux profile is assigned along the top of the domain and defined by:
\begin{align}\label{2DLaser}
	q_{2D}(x) &= .4\text{e}5\exp\left(\frac{-\left(L/2-x\right)^4}{1\text{e-}12}\right).
\end{align}
We fix the global mesh level $h_+=1/160$ and $\kappa_+=1.0$. 
\par The 3D is analogous, considering instead a cube of dimension $1/40 \times 1/40 \times 1/40$, and modifying the heat flux profile accordingly to: 
\begin{align}\label{3DLaser}
	q_{3D}(x,y) &= .4\text{e}5\exp\left(\frac{-\left(\left(L/2-x\right)^4 + \left(L/2-y\right)^4 \right) }{1\text{e-}12}\right).
\end{align}
\subsection{Experimental procedure}
Using the standardized problem setup defined previously, each experimental case is identified by the following characteristics:
\begin{itemize}
	\item Spatial dimension: $d$ = 2, 3;
	\item Degree of polynomial approximation $\mathcal{P}^m$: $m$= 1, 2;
	\item Local mesh resolution $h_-$: $1/320$, $1/640$, $1/1280$, $1/2560$. As global mesh resolution $h_+$ is fixed, these are identified by the ratio $h_+/h_-=$ $2$, $4$, $8$, $16$ respectively.
\end{itemize}
For brevity and clarity, we will refer to each case with this terminology: Case $d$=2, $m$=1, $h_+/h_-$=4 is understood as the two-dimensional case with linear polynomial approximation and $h_-$=$1/640$.
\par Each case then consists of 8 simulations, with $\kappa_-=.5^l$ for $l = 1,\,2,\,...,\,8$.\footnote{Due to memory constraints, the case $d$=3, $m$=2, $h_+/h_-$=16 was not considered.} We then compute linear and quadratic least-squares polynomial fittings between $\kappa_-/\kappa_+$ and $\specRad$: 
\begin{align}\begin{split}\label{linLeast}
\hat{\rho}_{L}\left(\frac{\kappa_-}{\kappa_+}\right) &\approx \left| a_1\left(\frac{\kappa_-}{\kappa_+}\right) +  a_0,\right|	\\
\hat{\rho}_{Q}\left(\frac{\kappa_-}{\kappa_+}\right) &\approx \left| b_2 \left(\frac{\kappa_-}{\kappa_+}\right)^2+b_1\left(\frac{\kappa_-}{\kappa_+}\right) +  b_0 \right|.	
\end{split}\end{align}
 If (\ref{scalingFn}) holds, we expect that $|b_2| \approx 0 $ and $a_1 \approx -a_0$, yielding:
 \begin{align}\label{linLeast2}
\hat{\rho}\left(\frac{\kappa_-}{\kappa_+}\right) &\approx \left|-a_1\left(\frac{\kappa_-}{\kappa_+}-1\right)\right|,
\end{align}
implying that:
\begin{align}\label{tildeLinLeast}
\widetilde{C} &\approx -a_1,
\end{align}
giving us an estimate of $\widetilde{C}$ for each case.
\par To summarize clearly, we will
\begin{itemize}
\item For each spatial dimension $d$, polynomial degree $m$, and local mesh resolution $h_+/h_-$, simulate the model problem with $\kappa_-=\frac{1}{2}^k$ for $k$=$1,\,2,\,...,\,8$. ;
\item The results of these simulations will be used to compute (\ref{linLeast}) and ideally estimate $\widetilde{C}$.
\end{itemize}
 
We start considering the unrelaxed case ($\theta=1$), for which we know from \eqref{cond30} and \eqref{scaling2}
that for \silvia{$\kappa_-/\kappa_+ < 2$} convergence \silvia{of the continuous scheme (\ref{eq:k0})-(\ref{eq:k2})} is guaranteed. Successively, we consider the impact of
the relaxation parameter $\theta$, where we find that the optimal relaxation parameter in the case of $\kappa_-/\kappa_+ \geq 2$ is given by $\kappa_+/\kappa_-$.

\subsection{Results}
\subsubsection{Convergence analysis}
We first seek to verify that $\specRad$ does indeed determine convergence, and that the scaling law (\ref{scalingFn}) is valid. We examine the case $d$=2, $m$=2, $h_+/h_-$=8, in detail, preferring a two-dimensional example in order to both run simulations and compute $\specRad$ quickly. 
\par
\begin{figure}
\centering
\includegraphics[scale=.25]{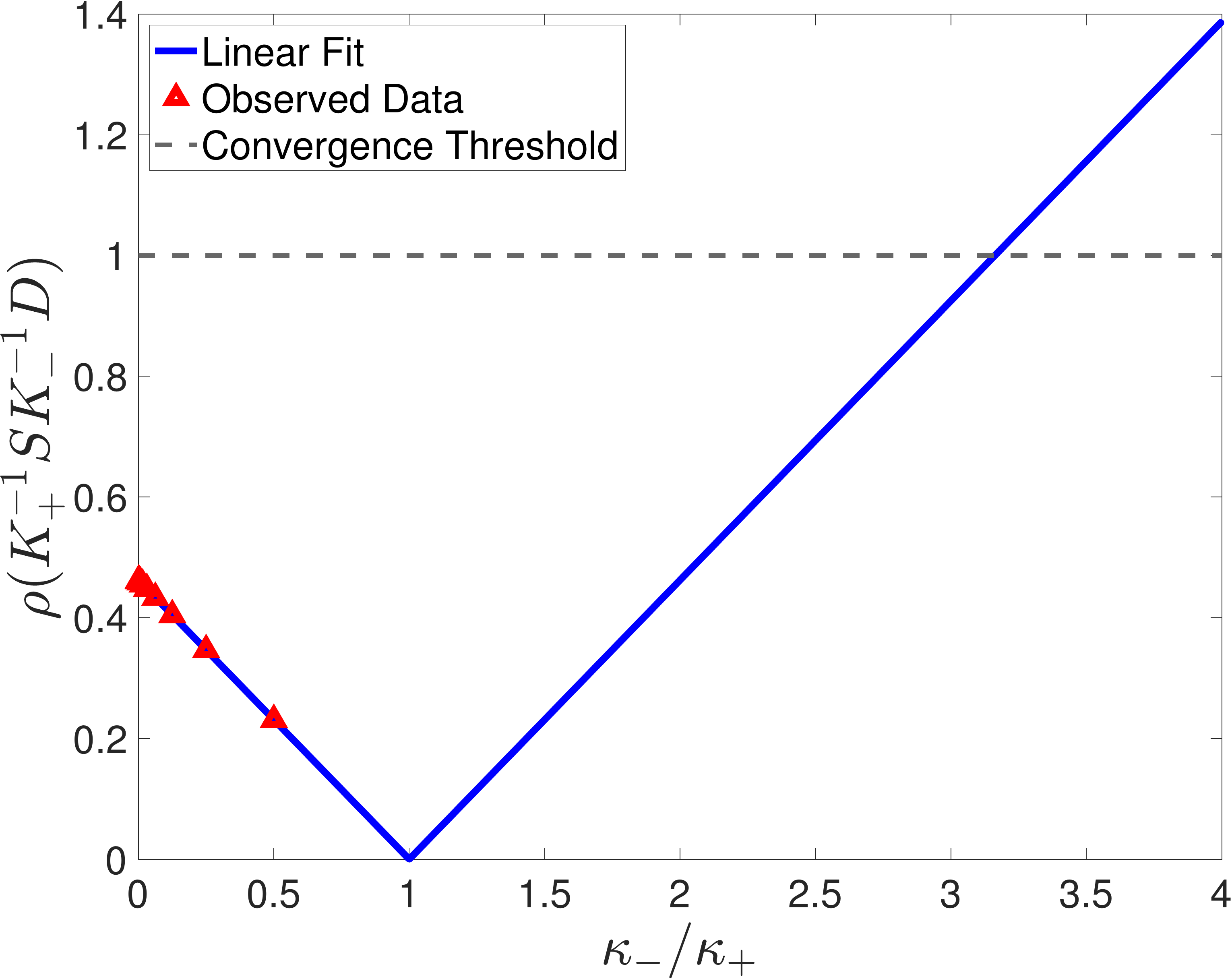}	\caption{Observed and predicted values of $\specRad$ for $d$=2, $m$=2, $h_+/h_-$=8. Divergence predicted for $\kappa_-/\kappa_+ \approx 3.16$.}\label{fig:eigFitting}
\end{figure}

 In Fig. \ref{fig:eigFitting}, we plot the computed $\specRad$ compared to $\kappa_-/\kappa_+$ from the simulation pool. The results of (\ref{linLeast}) give $a_1$=.4637, $a_0$=-.4637, and $b_2$=$8.84\text{e}-7$, confirming (\ref{scalingFn}) and implying that $\widetilde{C}\approx.4637$. 
\par As pointed out, we should observe convergence for all $\kappa_-/2 < \kappa_+$. The critical value for the divergence should be for $\widetilde{C} |\kappa_-/\kappa_+ - 1|
\approx 1$, that means in this case
 $\kappa_-/\kappa_+ = 1/\widetilde{C} + 1 \approx 3.16$. In Table \ref{tab:Data}, we report the number of iterations necessary for convergence, computed $\specRad$, and predicted $\specRad$ for a range of $\kappa_-/\kappa_+$. We define convergence as when relative difference in $L^2$ norm between consecutive iterations of $\boldsymbol{T}_+^{h,k}$ drops below 1e-8.
\par Referring to the data, we indeed observe the predicted behavior. The predicted and computed $\specRad$ are in perfect agreement, providing strong evidence for (\ref{scalingFn}). 
\par  This evidence strongly validates the theory detailed in previous sections. Note as well that, in general, $\kappa_- < \kappa_+$ for our preferred application of additive manufacturing, implying that for cases of this type convergence should not be a problem \cite{VBA2019}. 

\paragraph{Effects of relaxation}

\par
\begin{figure}
\centering
\includegraphics[scale=.25]{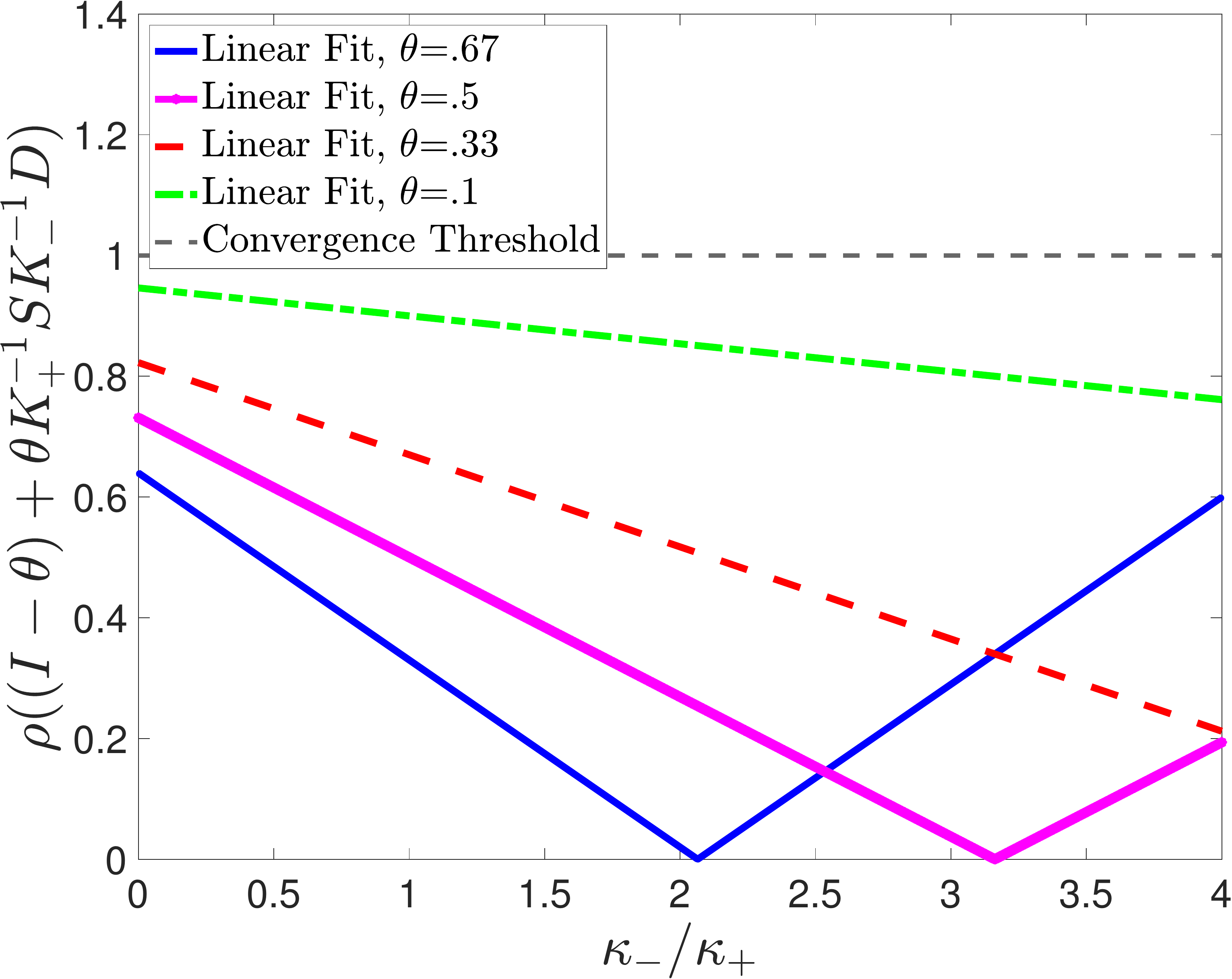}	\caption{Predicted values of $\rho\left( 	\left(1-\theta\right)I_+ + \theta K_+^{-1}SK_-^{-1}D \right)$ for $d$=2, $m$=2, $h_+/h_-$=8, various different values of $\theta$.}\label{fig:eigFitting2}
\end{figure}

 We note that the negative slope of the principal eigenvalue shown in Fig. \ref{fig:eigFitting} (note absolute value is shown, the signed value of the largest eigenvalue become negative) implies that the relaxation scheme (\ref{discConvCriterion}) will be effective for this problem. In Fig \ref{fig:eigFitting2}, we show the impact of $\theta$ on (\ref{discConvCriterion}). We therefore expect that, even for larger gaps for which Fig. \ref{fig:eigFitting} indicates divergence will occur,  in Fig. \ref{fig:eigFitting2}, appropriate $\theta$ can resolve this problem and ensuring convergence. We verify this numerically and, indeed, with $\theta$=.67 for $\kappa_-/\kappa_+$=3.16, the simulation converged in 12 iterations. We note also that Fig. \ref{fig:eigFitting2} suggests this parameter can be optimized depending on the problem, something that may be explored in future work.

\begin{table}
\centering
\begin{tabular}{|c|c|c|c|}
\hline
$\kappa_-/\kappa_+$ & Predicted $\rho$ & Actual $\rho$ & Num. iterations	\\ \hline
1.5 & .2319 & .2319 & 8 \\ \hline
2.0 & .4637 & .4637 & 14 \\ \hline
2.5 & .6956 & .6956 & 30 \\ \hline
3.0 & .9274 & .9274 & 143 \\ \hline
3.1 & .9738 & .9738 & 408 \\ \hline
3.15 & .997 & .997 & 3988 \\ \hline
3.16 & 1.001 & 1.001 & \textit{No convergence} \\ \hline
\end{tabular}\caption{Testing the relationship between $\kappa_-/\kappa_+$ and convergence. The predicted and computed $\specRad$ are in perfect agreement, and divergence occurs at the point predicted by $\specRad$. }\label{tab:Data}
\end{table}

\subsubsection{Spectral growth analysis}
\par Having confirmed with the previous that $\specRad$ determines the convergence and scales according to (\ref{scalingFn}), we now wish to analyze the influence of $d$, $m$, and $h_+/h_-$ on $\specRad$. As we have strong evidence that (\ref{scalingFn}) holds, we expect the influence of the polynomial degree, spatial dimension and mesh resolution on $\specRad$ to be reflected in the value of $\widetilde{C}$, predicted by (\ref{tildeLinLeast}). In particular, we are interested mostly in how this value grows as we refine $h_-$. Ideally we would like to observe that:
 \begin{align}\label{targetGrowth}
 \widetilde{C} &\sim \mathcal{O}\left(\log \frac{h_+}{h_-}\right	).
 \end{align}
In view of (\ref{scalingFn}), 
this indicates that as the local mesh is refined, the convergence dose not deteriorate, something of great practical importance for the application of the Two-level DD method.

\par
\begin{figure*}
        \centering
        \begin{subfigure}[b]{0.475\textwidth}
            \centering
            \includegraphics[width=\textwidth]{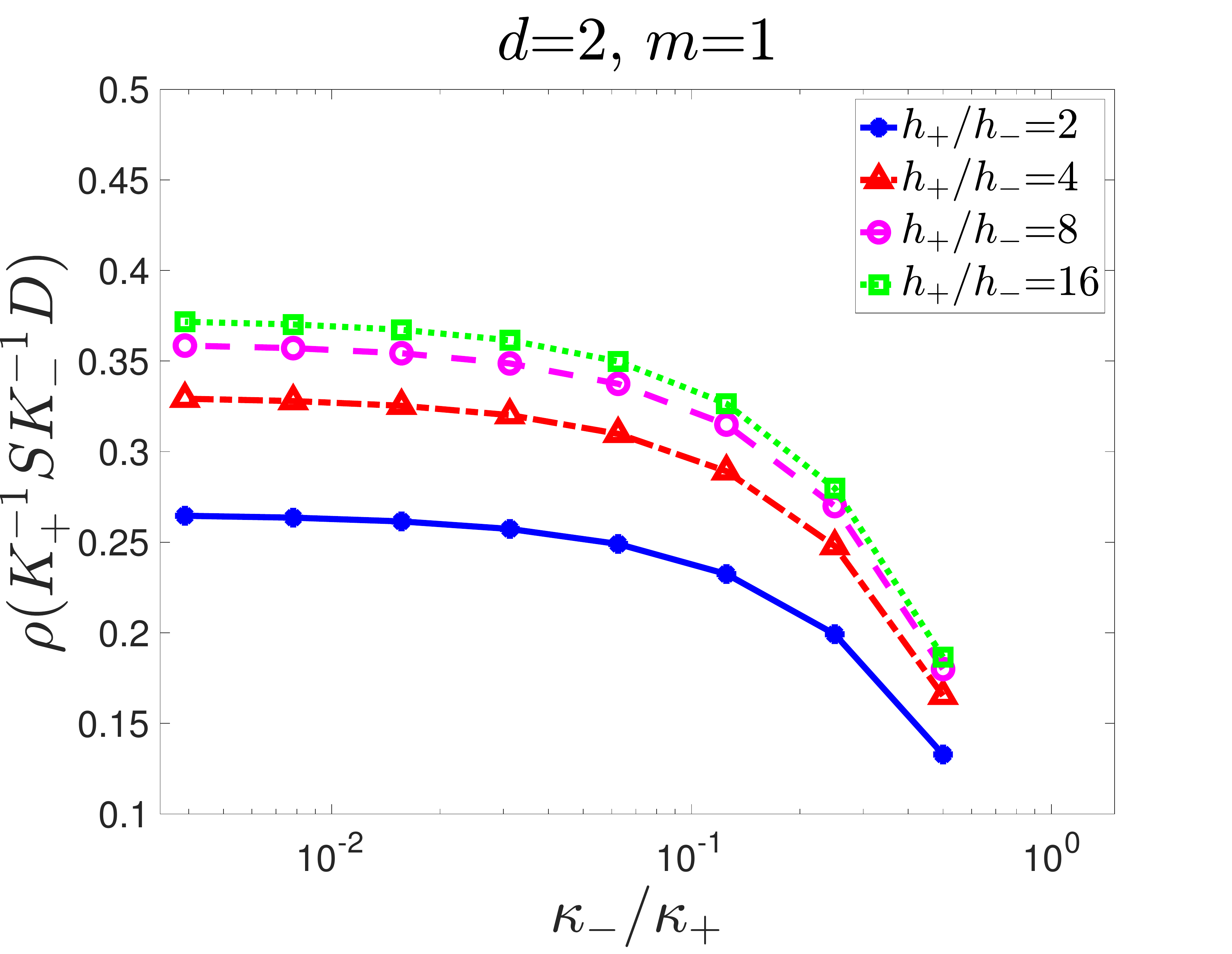}
  %          \caption[Network2]%
 %           {{\small Network 1}}    
%            \label{fig:mean and std of net14}
        \end{subfigure}
        \quad
        \begin{subfigure}[b]{0.475\textwidth}  
            \centering 
            \includegraphics[width=\textwidth]{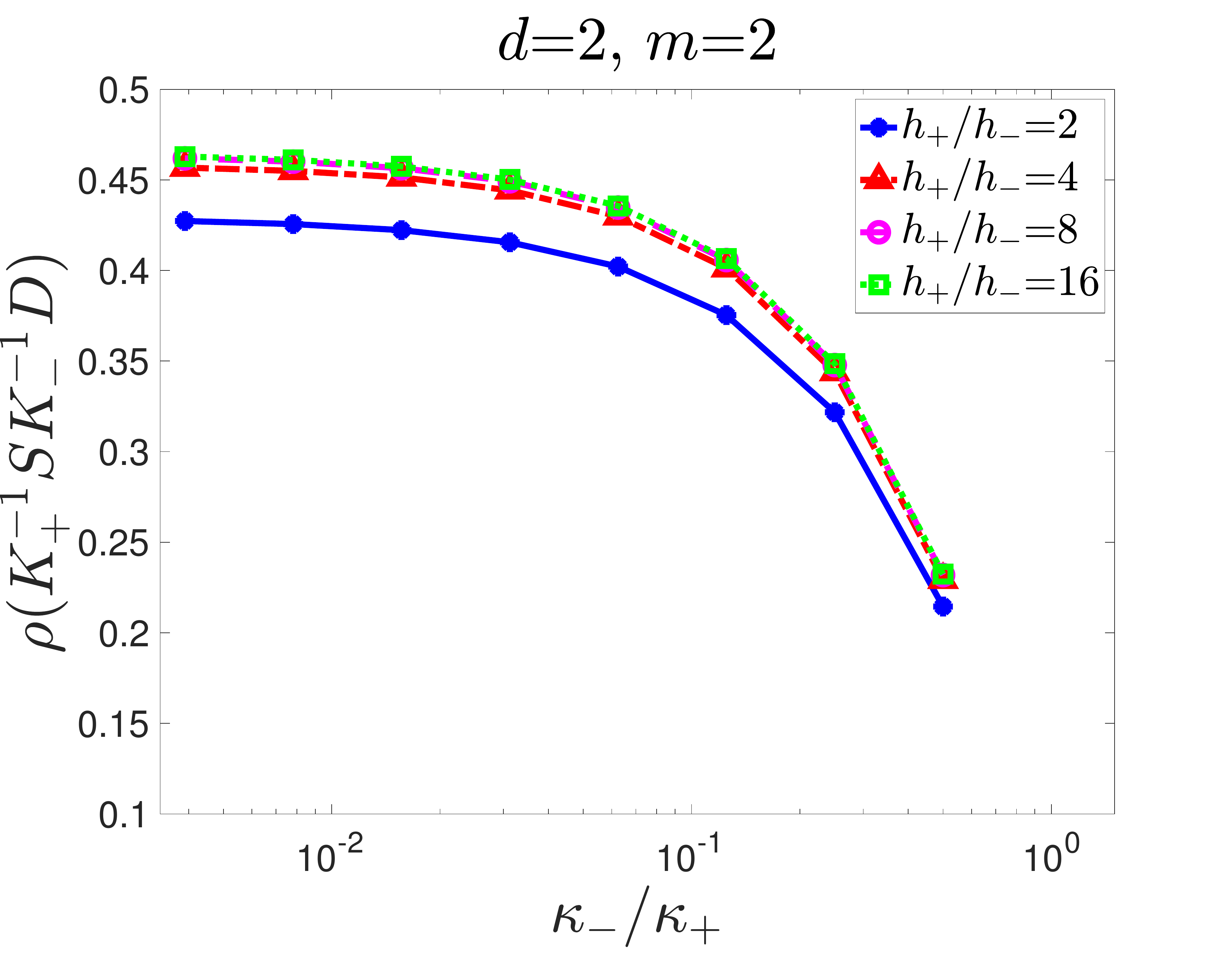}
   %         \caption[]%
  %          {{\small Network 2}}    
 %           \label{fig:mean and std of net24}
        \end{subfigure}
        \vskip\baselineskip
        \begin{subfigure}[b]{0.475\textwidth}   
            \centering 
            \includegraphics[width=\textwidth]{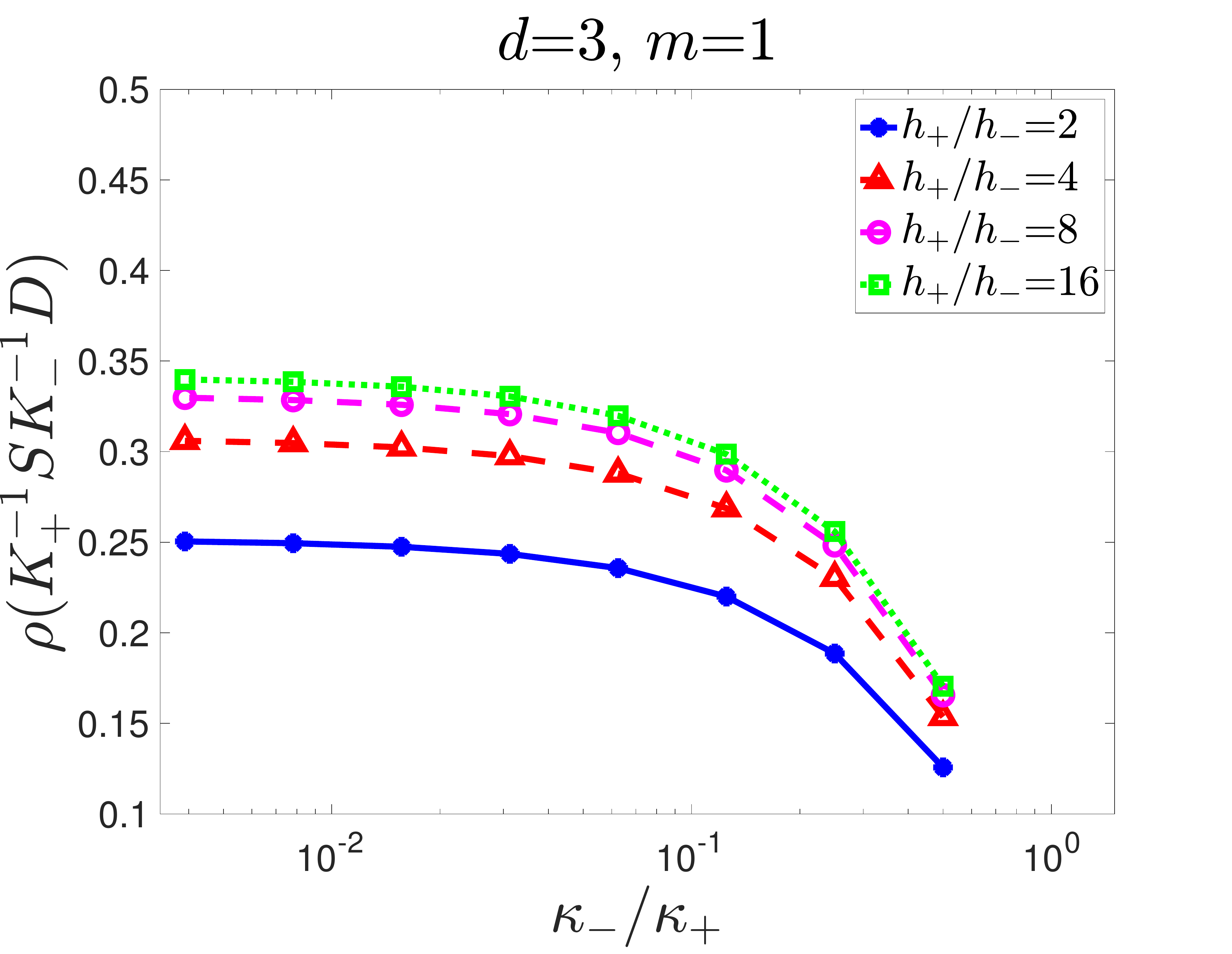}
  %          \caption[]%
 %           {{\small Network 3}}    
%            \label{fig:mean and std of net34}
        \end{subfigure}
        \quad
        \begin{subfigure}[b]{0.475\textwidth}   
            \centering 
            \includegraphics[width=\textwidth]{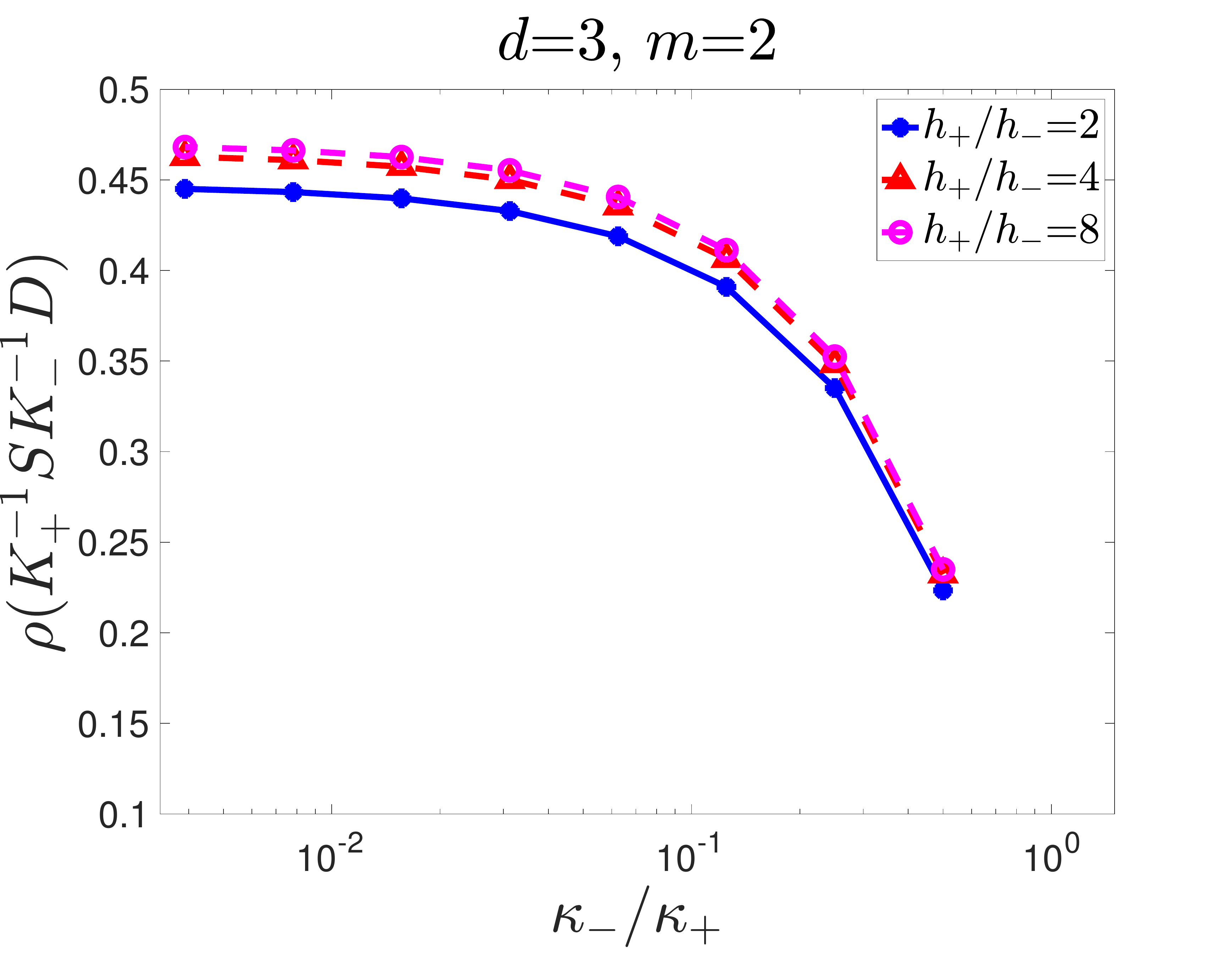}
   %         \caption[]%
  %          {{\small Network 4}}    
 %           \label{fig:mean and std of net44}
        \end{subfigure}
        \caption[ Relationship between $\kappa_-/\kappa_+$ for different $d$, $m$ and $h_+/h_-$. The top row corresponds to 2D cases, the bottom to 3D; the left column corresponds to linear polynomial approximations, the right to quadratic. We observe logarithmic increase in $\specRad$ for each $d$ and $m$.  ]
  			{\small  Relationship between $\kappa_-/\kappa_+$ for different $d$, $m$ and $h_+/h_-$. The top row corresponds to 2D cases, the bottom to 3D; the left row corresponds to linear polynomial approximations, the right to quadratic. For each $d$ and $m$, the increase in $\specRad$ with $h_+/h_-$ appears logarithmic. }
        \label{fig:growthPlots}
    \end{figure*}

 In Fig. \ref{fig:growthPlots} we display the relationship between $\specRad$ over the different $d$, $m$, and $h_+/h_-$. The plots appear indicate that, for a given $d$ and $m$, $\specRad$ displays a logarithmic growth trend, as hoped.  
    
    \begin{figure}
\centering
\includegraphics[scale=.35]{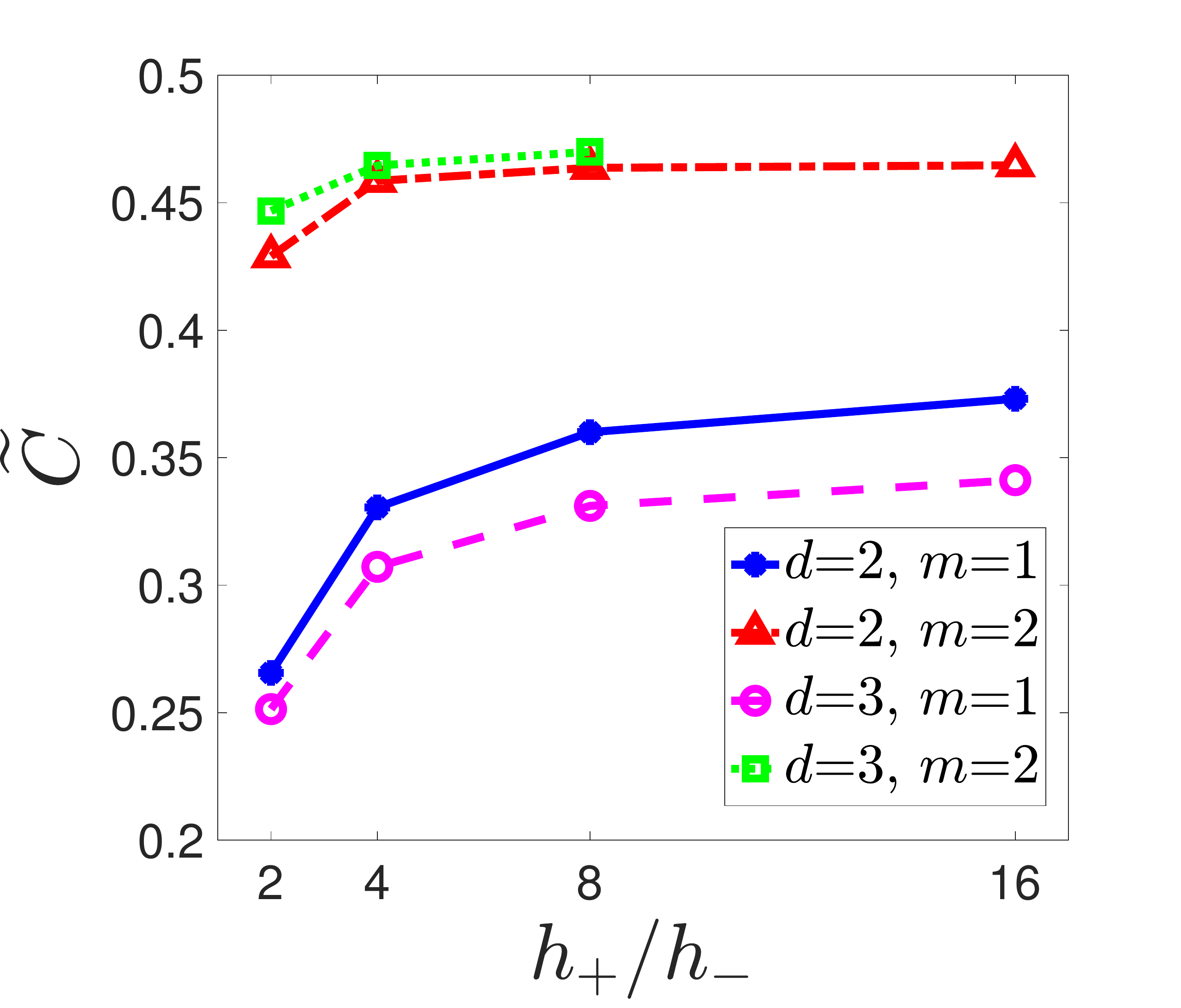}	\caption{Growth of $\widetilde{C}$ as $h_-$ is refined across different $d$ and $m$. For all cases, we observe a logarithmic growth trend, suggesting that the performance of the algorithm is robust to refinement of the local mesh. }\label{fig:ctildeplot}
\end{figure} We confirm this logarithmic scaling behavior further in Fig. \ref{fig:ctildeplot}, which plots the computed value of $\widetilde{C}$ against different $h_+/h_-$ levels for each $d$ and $m$. The evidence here suggests that the degree of polynomial approximation influences the spectral behavior far more than spatial dimension, at least for this problem. In particular, we note that when $m$=2, the value of $\widetilde{C}$ is larger than when $m$=1, yet it grows very little.    
\par The numerical experiments performed suggest several conclusions. First and foremost, we observe strong validation of our theoretical prediction; namely, that the convergence of the Two-level DD method is determined by $\rho\left(K_+^{-1}SK_-^{-1}D\right)$. Further, this quantity depends primarily on $\kappa_+$ and $\kappa_-$. 
Trivially,
for $\kappa_+ \approx \kappa_-$ the convergence is fast, as this cases approaches the no-jump case in the coefficients, perfectly approximated by a monolithic approach. Looking over the entire dataset, we observed the dependence of $\specRad$ on the local mesh resolution to be logarithmic, a critical observation, as we do not expect rapid performance deterioration with mesh refinement.  
\subsection{Comparison with monolithic solver on a fitted mesh}
\par \alex{In this section, we now compare the performance of the Two-level DD method with the use of a monolithic solver in the entire domain. One of the clear advantages of the Two-level DD method is in its ability to employ two uniform meshes, each with homogeneous values of $\kappa$, potentially leading to better conditioning properties. Such improvements in conditioning, when sufficiently large, may render the use of a DD method competitive, even despite the need to iterate. }
\par \alex{We will proceed by considering the same basic two-dimensional setup shown in the previous sections, now using Jacobi-preconditioned GMRES for all linear solves. For the purposes of comparison, we consider the solution time and necessary GMRES iterations for the Two-level solves and for the monolithic solution over a range of $\kappa_-/\kappa_+$ and $h_+/h_-$, fixing $\kappa_+$ as 1.0,$m=2$. The monolithic solution is performed on a fitted mesh, such that the mesh features a size of $h_-$ in $\Omega_-$ and $h_+$ in $\Omega_+$, where all nodes along the boundary between $\Omega_-$ and $\Omega_+$ are conformal to the separation in the physical domain  }

\begin{table}
\tiny
\centering
\begin{tabular}{|c|c|c|c|c|c|c|}
\hline
$\kappa_-/\kappa_+$ & $h_+/h_-$  &Local GMRES Its. &Global GMRES Its. & Monolithic GMRES Its.  & Time (Two-lvl) & Time (monolithic)	\\ \hline
1 & 8 & 73 & 290 & 3608 & .69s & 7.9s \\ \hline
.1 & 8 & 77  & 310 & 1557 & 1.11s & 3.56s \\ \hline
.01 & 8 & 99 & 382 & 1854 & 1.34s & 4.16s \\ \hline
1 & 16 & 176 & 290 & 6302 & 2.88s & 26.72s \\ \hline
.1 & 16 & 137  & 309 & 1916 & 4.61s & 8.59s \\ \hline
.01 & 16 & 222 & 381 & No convergence & 5.74s & NA \\ \hline
\end{tabular}\caption{Comparison with monolithic solver: observing the relationship between $\kappa_-/\kappa_+$, $h_+/h_-$, and numerical solution behavior in terms of linear solve difficulty and solution time.  }\label{tab:monolithicSolve}
\end{table}
\par \alex{We show the results in Table \ref{tab:monolithicSolve}. We generally observe faster solution times and better behavior with respect to changes in $\kappa_-/\kappa_+$ for the Two-level DD solution. Notably, for $h_+/h_-$=16 and $\kappa_-/\kappa_+$=.01, the Jacobi-preconditioned GMRES fails to converge for the monolithic solution, while converging in a reasonable amount of iterations for the Two-level DD case. The behavior with respect to mesh size is a general increase in necessary iterations for each method, as expected. However, for $\kappa_-/\kappa_+$ the effect is less clear; it is clear, however, that the Two-level DD method appears only mildly sensitive in this respect, with the number of necessary GMRES iterations remaining around the same order of magnitude in all cases. In contrast, the monolithic method shows highly variable and unpredictable behavior, likely a result of conditioning issues brought on by the large heterogeneities present in the problem.
}
\subsection{Extension to the nonlinear case}
\par \alex{We will now consider the same general problem setup, but will introduce nonlinear problems in which the thermal conductivity depends on the temperature, with the nature of this dependence different in the local and global domains. In \cite{VBA2019}, it was shown that, through some modification of the transmission terms, the Two-level method is also consistent with the original problem formulation in such a case. We refer the reader to this work for the additional details. We note in this instance the necessity to iterate no longer represents a serious issue, as the nonlinear nature of the problem requires that some sort of iterative method be employed, regardless of whether one solves in a monolithic fashion or using a Two-level DD method. While this test is, strictly speaking, outside the analysis shown within the present work, we nonetheless feel its inclusion is important, given the eventual desired application of the proposed method to more complex problems.}
\par \alex{In terms of setup, the problem is identical to those shown in the previous sections, with the important distinction now that the thermal conductivity in the local domain is considered as the powder form of stainless steel 316L (see e.g. \cite{VA2019, Mills2002}). For the global problem, we define $\kappa_+$ as the solid form of stainless steel 316L in $\Omega_A$. This now introduces an additional challenge, as we must extend the definition of $\kappa_+$ in $\Omega_B$ in some way. We denote this extension as $\kappa_{+,B}$. For ease of computation, it would be convenient to consider $\kappa_{+,B}$ as constant, and we will examine the viability of this choice over a range of different constants. }
\par \alex{In each case, we solve a nonlinear problem using a Picard-type linearization and each linear system using GMRES with a Jacobi-style preconditioner. We compare the Two-level DD solution to a monolithic solution on a conformal nonuniform mesh in terms of GMRES iterations, nonlinear iterations, and overall solution time.}

\begin{table}
\centering
\begin{tabular}{|c|c|c|c|}
\hline
$\kappa_{+,B}$ & Avg. Global GMRES Iterations & Num. nonlinear iterations  & Time (s)	\\ \hline
.1 & 643 &  13 & 26.8 \\ \hline
.2 & 621 & 11 & 26.1 \\ \hline
.3 & 707 & 10 & 25.6 \\ \hline
.4 & 728 & 12 & 28.15 \\ \hline
.5 & 725 & 11 & 28.01 \\ \hline
.6 & 721 & 12 & 28.87 \\ \hline
.7 & 779 & 14 & 30.42 \\ \hline
.8 & 875 & 15 & 32.21 \\ \hline

\end{tabular}\caption{Nonlinear problem, observing the relationship between the constant extension $\kappa_{+,B}$ of $\kappa_B$ and the difficulty of the linear and nonlinear solves. We see that increasing $\kappa_{+,B}$ generally leads to more difficult linear solves, while a $\kappa_{+,B}$ closer to .3 minimizes the number of necessary nonlinear iterations. }\label{tab:NonlinearData}
\end{table}

\par \alex{In Tab. \ref{tab:NonlinearData}, we display the results of these tests. Solving the monolithic reference problem, we note that the mean value of $\kappa_+$ in near $\Omega_B$ is approximately .306, and hence we expect values for $\kappa_{+,B}$ close to .306 to provide superior numerical performance.  In accordance with our expectation, we indeed find that extending $\kappa_{+,B}$ as a constant such that the jump between $\kappa_{+,A\setminus B}$ and $\kappa_{+,B}$ is small reduces the necessary number of nonlinear iterations for convergence. As for the difficulty of linear solves, the behavior is less clear; lower values of $\kappa_{+,B}$ seem to result in easier linear solves, though not necessarily fewer nonlinear iterations. For the purposes of comparison, the monolithic reference solution required 8 nonlinear iterations, each such iteration requiring an average of 6354 GMRES iterations, resulting in an overall solve time of 126.3 seconds. In this test, we observe both superior conditioning and solver performance when compared to the monolithic approach, as well as an important numerical confirmation that extending $\kappa_B$ as a constant inside $\Omega_B$ may be a reasonable choice. Based on the results of these tests, if one does not have a good estimate for what constant value to use, in general lower values seem to provide better linear and nonlinear conditioning behavior.  }

\section{Conclusions and Future Work}\label{sec:concl}
With this work, \silvia{we present a  theoretical analysis of the Two Level DD method for the solution of heterogeneous material problems, obtained by combining a non-overlapping 
Dirichlet-Neumann DD technique with a fictitious domain approach.} Under the simplifying assumption that the coefficients are piece-wise constants, we have established that, at the continuous level, the considered formulation converges to the desired solution. At the discrete level, we showed that one may interpret an algorithm of this type as a block Gauss-Seidel type iteration or as a truncated Neumann series, from which we were able to postulate a convergence criterion and its scaling behavior. We then performed several two- and three-dimensional simulations which validated our predictions. Importantly, \silvia{our numerical tests show} that the convergence of the method is expected to be robust with respect to local refinement, which is important information of great practical interest for the application of such algorithms. To demonstrate the potential of the approach, we also performed some test in the nonlinear case, and made two important conclusions: that the superior conditioning afforded by the DD method provides a numerical advantage over applying a fitted non uniform mesh with a monolithic approach, and that one may extend $\kappa_+$ as a constant inside $\Omega_B$ while maintaining good numerical performance, provided such a value is chosen carefully.
\par 
%Though this work has greatly extended the theory of the Dirichlet-Neumann domain decomposition methods in heterogeneous media,
\silvia{While this work is a first step in understanding the theoretical and numerical features of the approach considered, many important questions remain. While we numerically examined the behavior of the method on nonlinear problem, in the theoretical analysis we only consider the steady problem and we rely on the simplifying assumption that the coefficients are piece-wise constants. The techniques shown here should extend promptly to unsteady problems with constant coefficients, however, for non-constant coefficients and for nonlinear problems, additional care must be taken in both the design of the method and convergence analysis.} This extends to the discrete problem as well, as the Gauss-Seidel/Neumann series interpretations shown here will require significant adaptation. \silvia{Finally, the approach here considered can be extended in different directions. To achieve higher efficiency, other domain decomposition strategies may be used in combination with fictitious domains, instead of the simple Dirichlet-Neumann iterations here considered. To achieve optimal error, one might consider the possibility of giving up the requirement that $T^+ = T$ on $\Omega_-$ and look for an extension $f$ for which $T^+$ is smoother than $T$ (which has a jump in the normal derivative along $\gamma_i$, thus limiting the $H^1(\Omega_+)$ error to $h^{1/2}$).}

\section{Acknowledgments}
\textit{This work was partially supported by the Italian Minister of University and Research through the project "A BRIDGE TO THE FUTURE: Computational methods, innovative applications, experimental validations of new materials and technologies” (No. 2017L7X3CS) within the PRIN 2017 program, as well as Regione Lombardia through the project "MADE4LO - Metal ADditivE for LOmbardy" (No. 240963) within the POR FESR 2014-2020 program.}

\bibliographystyle{abbrv}  \bibliography{Viguerie_Veneziani.bib}
\end{document}